\documentclass[reqno,centertags,12pt,a4paper]{amsart}
\usepackage{dsfont}
\usepackage{amscd}
\usepackage{amsmath,tikz}
\usepackage{amssymb}
\usepackage{amsfonts}
\usepackage{hyperref}
\usepackage{enumerate}
\usepackage{latexsym}
\usepackage{cases}
\usepackage{amsthm}
\usepackage{graphicx}
\usepackage{verbatim}
\usepackage{mathrsfs}
\usepackage{color}
\usepackage{wasysym}
\usepackage{makecell}
\newcommand\pxx{%
	\mathrel{\text{\tikz[baseline] \draw (0em,-0.3ex) -- (.4em,1.7ex) (.2em,-0.3ex) -- (.6em,1.7ex);}%
}}

\newtheorem{lemma}{Lemma}[section]
\newtheorem{proposition}[lemma]{Proposition}
\newtheorem{corollary}[lemma]{Corollary}
\newtheorem{definition}[lemma]{Definition}

\newtheorem{remark}[lemma]{Remark}
\newtheorem{theorem}[lemma]{Theorem}
\theoremstyle{definition}
\numberwithin{equation}{section}

\allowdisplaybreaks[2] 

\begin{document}

\title[Fourier decay of fractal measures]{Fourier decay of fractal measures on surfaces of co-dimension two in $\mathbb{R}^5$}
\author[Zhenbin~Cao,~Changxing ~Miao~and~Zijian~Wang]{Zhenbin~Cao,~Changxing ~Miao~and~Zijian~Wang}

\date{\today}

\address{Institute of Mathematics, Henan Academy of Sciences, Zhengzhou 450046, China}
\email{caozhenbin@hnas.ac.cn}

\address{Institute of Applied Physics and Computational Mathematics, Beijing 100088, China}
\email{miao\_{}changxing@iapcm.ac.cn}

\address{ Zhongkai University of Agriculture and Engineering, Guangzhou 510225, China   }
\email{ saswangzijiana@126.com    }


\begin{abstract}
In this paper, we study the Fourier decay of fractal measures on the quadratic surfaces of high co-dimensions. Unlike the case of co-dimension 1, quadratic surfaces of high co-dimensions possess some special scaling structures and degenerate characteristics. We will adopt the strategy from Du and Zhang \cite{DZ}, combined with the broad-narrow analysis with different dimensions as divisions, to obtain a few lower bounds of Fourier decay of fractal measures on quadratic surfaces of co-dimension two in $\mathbb{R}^5$.

\end{abstract}

\maketitle

\section{INTRODUCTION}\label{section1}

Let $m\geq 2$. For any Borel set $A\subset\mathbb{R}^m$, we denote by $\mathcal{M}(A)$ the collection of all Borel finite measures with compact support contained in $A$. 
For $\alpha\in(0,m]$, we define the $\alpha$-dimensional energy of $\mu\in\mathcal{M}(B^m(0,1))$ as
$$I_{\alpha}(\mu):=\int\int\frac{d\mu(x)d\mu(y)}{|x-y|^{\alpha}}=C_{\alpha,m}\int|\xi|^{\alpha-m}|\widehat{\mu}(\xi)|^2d\xi,$$
where $\widehat{\mu}$ is defined by 
$$\widehat{\mu}(\xi)=\int e^{- 2\pi ix\cdot \xi}d\mu(x).$$
In general, for $\mu\in\mathcal{M}(B^m(0,1))$, $I_{\alpha}(\mu)<\infty$ does not imply  
$$|\widehat{\mu}(\xi)|\rightarrow0,\quad \quad\text{as~} |\xi|\rightarrow\infty.$$ 
However,  the integral average of $\widehat{\mu}(R\xi)$ on the compact surface $\Pi$ tends to zero as $R \rightarrow \infty$, i.e. 
\begin{align}\label{12121}
	\int_{\Pi}|\widehat{\mu}(R\xi)|^2d\sigma(\xi)\rightarrow 0, \quad\quad\text{as~} R\rightarrow \infty,
\end{align}
where $d\sigma$ is the surface measure on $\Pi$. Motivated by the research on geometric measure theory and partial differential equations, more and more people are interested in the question how fast decay does the average of $\widehat{\mu}(R\xi)$ over different compact surfaces have.
Let $\beta(\alpha,\Pi)$ denote the supremum of the number $\beta>0$ for which
\begin{align}\label{WZJ1}
	\int_{\Pi}|\widehat{\mu}(R\xi)|^2d\sigma(\xi)\lesssim R^{-\beta},\quad\quad\forall ~R>1,
\end{align}
for all $\mu\in\mathcal{M}(B^m(0,1))$ with $I_{\alpha}(\mu)<\infty$. 

A major way to finding the exact values of  $\beta(\alpha,\Pi)$ is to build the weighted restriction estimate associated with $\Pi$. Let $d,n\geq 1$ with $d+n=m$, and $\phi: [0,1]^d \rightarrow \mathbb{R}^n$ be the parameter representation of $\Pi$. 
Define the Fourier extension operator
$$E^{\Pi}f(x,y):=\int_{[0,1]^d} e^{ 2\pi i(x\cdot\xi+y\cdot\phi(\xi))}f(\xi)d\xi,\quad \quad (x,y)\in \mathbb{R}^d\times \mathbb{R}^n.$$
A measure $\mu\in\mathcal{M}(B^{m}(0,1))$ is called $\alpha$-dimensional if it satisfies
\begin{align*}
	\mu(B(b,r))\leq C_\mu r^{\alpha},\quad\quad\forall ~b\in\mathbb{R}^{m},~~\forall ~r>0.
\end{align*}
Denote  $\mu_R(\cdot):=R^{\alpha}\mu(R^{-1}\cdot)$. Define $s(\alpha,\Pi)$ be the smallest constant $s$ such that the inequality
\begin{align}\label{qwe}
	\|E^{\Pi}f\|_{L^2(d\mu_R)}\lesssim R^{s}\|f\|_{L^2}, \quad \quad \forall ~R \geq 1,
\end{align}
holds for all $\alpha$-dimensional measures $\mu\in\mathcal{M}(B^{m}(0,1))$ and every $f\in L^2(\mathbb{R}^d)$. These two questions are equivalent \cite{BEH20} in the following sense:
\begin{align}\label{W1}
	\beta(\alpha,\Pi)=\alpha-2s(\alpha,\Pi).
\end{align}

We give one brief introduction about the research history of $	\beta(\alpha,\Pi)$. The most important case is $\Pi=\mathbb{S}^{m-1}$, which is connected with Falconer's distance set conjecture \cite{F85} in geometric measure theory. In fact, if (\ref{WZJ1}) holds for $\beta=m-\alpha$, then it implies that the distance set of each compact set whose Hausdorff dimension is larger than $\alpha$ has positive Lebesgue measure. When $m=2$, 
Mattila \cite{M87} obtained the sharp bound of $\beta(\alpha,\mathbb{S}^{1})$ when $\alpha \in (0,1]$ by the asymptotic expression of the Bessel function. His argument is also fitted for general $m \geq 3$ when $\alpha\in(0,m/2]$ (see (\ref{his1})). Wolff \cite{W99} built the corresponding weighted restriction estimate by making use of the Kakeya estimate and square function estimate in the plane, which implies the sharp bound of $\beta(\alpha,\mathbb{S}^1)$ when $\alpha\in(1,2]$. The exact values of $\beta(\alpha,\mathbb{S}^{1})$ are
\begin{equation}
	~\beta(\alpha,\mathbb{S}^{1})= \begin{cases}
		~	\alpha, \quad \quad&\alpha\in(0,\frac{1}{2}],\quad\quad $(Mattila \cite{M87})$   \\
		~	\frac{1}{2} , \quad \quad& \alpha\in(\frac{1}{2},1],\quad\quad $(Mattila \cite{M87})$  \\
		~	\frac{\alpha}{2}, \quad  \quad&  \alpha\in(1,2],\quad\quad \,$(Wolff \cite{W99})$.
	\end{cases}
\end{equation}
For the counterexamples, readers can see \cite{M87,W99}. When $m \geq 3$, Du and Zhang \cite{DZ} established the fractal $L^2$ restriction estimate, which implies the best lower bound of $\beta(\alpha,\mathbb{S}^{m-1})$ so far when $\alpha\in (m/2,m]$, see also \cite{DGOWWZ20} when  $m=3$. The current best lower bounds of $\beta(\alpha,\mathbb{S}^{m-1})$ are
\begin{equation}\label{his1}
	~\beta(\alpha,\mathbb{S}^{m-1})\geq \begin{cases}
		~	\alpha,  \quad&\alpha\in(0,\frac{m-1}{2}],\quad\,\,$(Mattila \cite{M87})$   \\
		~	\frac{m-1}{2} , \quad& \alpha\in(\frac{m-1}{2},\frac{m}{2}], \quad$(Mattila \cite{M87})$ \\
		~	\frac{(m-1)\alpha}{m},  \quad&  \alpha\in(\frac{m}{2},m],\,\,\,\,\,\quad$(Du and Zhang \cite{DZ}$).
	\end{cases}
\end{equation}
The case when  $\Pi$ is the truncated light cone was studied in \cite{E04,CHL17,H19}.    Recently, Barron, Erdo\v{g}an and Harris \cite{BEH20} studied the case when $\Pi$ is the   hyperboloid.

In this paper we consider the case when $\Pi$ are the surfaces of high co-dimensions. The study on the surfaces of high co-dimensions is taken more and more seriously, since they have important applications in many areas. For example, the decoupling theory for the moment curve is connected with Vinogradov's Mean Value Theorem \cite{Woo16,BDG16} in number theory; the decoupling theory for the manifolds of co-dimension 2 in $\mathbb{R}^4$ is connected with Lindel$\ddot{\text{o}}$f hypothesis \cite{B17} in number theory; The Fourier restriction estimates and weighted Fourier restriction estimates for  general manifolds of co-dimension 2 are useful for obtaining Carleman's inequalities and weighted Carleman's inequalities \cite{De00} in PDEs, respectively.

Let $d,n \geq 1$. We denote by $\textbf{Q}(\xi)=(Q_1(\xi),...,Q_n(\xi))$ an $n$-tuple of real quadratic forms in $d$ variables. The graph of such a tuple, $S_{\textbf{Q}}=\{ (\xi,\textbf{Q}(\xi))\in [0,1]^d \times \mathbb{R}^n \}$, is a $d$-dimensional submanifold of $\mathbb{R}^{d+n}$. There have been a number of works studying the problems in harmonic analysis for the quadratic surfaces of high co-dimensions. Readers can see \cite{C85,C82,M96,BL04,BLL17,LL20,GO} for the restriction problems, and \cite{BD16,DGS19,GZK20,GORYZK19,O18} for the decoupling theory. No matter which problem, one  of the difficulties in the setting of co-dimension bigger than one is the lack of a notion of ``curvature". Recently, Guo, Oh, Zhang and Zorin-Kranich \cite{GOZZK} introduced one quantity $\mathfrak{d}_{d',n'}(\textbf{Q})$ (defined below), which can be used to measure the properties of ``curvature" on such surfaces, and obtained the corresponding sharp decoupling inequalities. 

For a tuple $\textbf{Q}(\xi)=(Q_1(\xi),...,Q_n(\xi))$ of quadratic forms with $\xi \in \mathbb{R}^d$, we define
$$  \text{NV}(\textbf{Q}) := \# \{    1 \leq d' \leq d : \partial_{\xi_{d'}} Q_{n'} \not\equiv 0 \text{~for~some~} 1\leq n'\leq n \}      .  $$
This quantity refers to ``the number of variables on which $\textbf{Q}$ depends". For $0 \leq d' \leq d$ and $0\leq n'\leq n$, define
\begin{equation}\label{main sign}
 \mathfrak{d}_{d',n'}(\textbf{Q}):= \inf_{\substack{   M \in \mathbb{R}^{d\times d} \\ \text{rank}(M)=d'   }}  \inf_{\substack{   M' \in \mathbb{R}^{n'\times n} \\ \text{rank}(M')=n'   }} \text{NV}(M'\cdot (\textbf{Q}\circ M)) ,
\end{equation}
where $\textbf{Q}\circ M$ is the composition of $\textbf{Q}$ with $M$. This definition can help us to give the classification of quadratic forms. We say that $\textbf{Q}$ and $\textbf{Q}'$ are equivalent, and denote by $\textbf{Q}\equiv\textbf{Q}'$ if there exist two invertible real matrices $M_1 \in \mathbb{R}^{d\times d}$ and $M_2 \in \mathbb{R}^{n\times n}$ such that 
$$   \textbf{Q}'(\xi) =  \textbf{Q}(M_1 \cdot \xi) \cdot M_2, \quad \quad \forall~ \xi \in \mathbb{R}^d. $$
A simple fact is that $\textbf{Q}\equiv\textbf{Q}'$ implies $ \mathfrak{d}_{d',n'}(\textbf{Q})= \mathfrak{d}_{d',n'}(\textbf{Q}')$ for all  $0 \leq d' \leq d$ and $0\leq n'\leq n$. We say that $\textbf{Q}$ satisfies the (CM) condition if
\begin{equation}
\int_{\mathbb{S}^{n-1}} |\text{det}(y_1 Q_1+...+y_n Q_n)|^{-\gamma}d\sigma(y) <\infty, \quad \quad \forall~0<\gamma<\frac{n}{d},
\end{equation} 
where each $Q_j$, $j=1,...,n$, denotes the matrix associated with the quadratic form $Q_j(\xi)$.
This condition was introduced by Christ \cite{C82,C85} and Mockenhaupt \cite{M96}, where they obtained the sharp Stein-Tomas type's inequalities for surfaces satisfying the (CM)  condition when $n=2$.

When $d=3$ and $n=2$, Guo and Oh classified quadratic forms according to some geometric conditions in \cite[Theorem 1.2]{GO}. We can rewrite their result by making use of the notations of \cite{GOZZK} as follows. The equivalence of the two theorems will be showed in Appendix 2.
\begin{theorem}\label{th1}
Let $d=3$ and $n=2$. Suppose that the quadratic form $\textbf{Q}=(Q_1,Q_2)$ satisfies   $\mathfrak{d}_{3,2}(\textbf{Q})=3$. Then 
\begin{itemize}
	\item[(a)] If $\textbf{Q}$ satisfies $\mathfrak{d}_{3,1}(\textbf{Q})=1$ and $\mathfrak{d}_{2,2}(\textbf{Q})=1$, then $\textbf{Q}\equiv (\xi_1^2,\xi_2^2+\xi_1\xi_3)$ or  $(\xi_1^2,\xi_2\xi_3)$.
	\item[(b)] If $\textbf{Q}$ satisfies $\mathfrak{d}_{3,1}(\textbf{Q})=1$ and $\mathfrak{d}_{2,2}(\textbf{Q})=2$, then $\textbf{Q}\equiv (\xi_1^2,\xi_2^2+\xi_3^2)$.
	\item[(c)] If $\textbf{Q}$ satisfies $\mathfrak{d}_{3,1}(\textbf{Q})=2$ and $\mathfrak{d}_{2,2}(\textbf{Q})=0$, then $\textbf{Q}\equiv (\xi_1 \xi_2,\xi_1\xi_3)$.
	\item[(d)] If $\textbf{Q}$ satisfies $\mathfrak{d}_{3,1}(\textbf{Q})=2$ and $\mathfrak{d}_{2,2}(\textbf{Q})=1$, then $\textbf{Q}\equiv (\xi_1\xi_2,\xi_2^2+\xi_1\xi_3)$ or $ (\xi_1\xi_2,\xi_1^2\pm\xi^2_3)$.
	\item[(e)] $\textbf{Q}$ satisfies $\mathfrak{d}_{3,1}(\textbf{Q})=2$ and $\mathfrak{d}_{2,2}(\textbf{Q})=2$ if and only if $\textbf{Q}$ satisfies the ${\rm (CM)}$ condition. 
\end{itemize}
\end{theorem}


Now we state the main Fourier decay results of the paper.

\begin{theorem}\label{th2}
	Let $\alpha\in(0,5]$. Then
	
	\begin{itemize}
		\item[(1)] For $\textbf{Q}$ satisfying ${\rm(a)}$ in Theorem \ref{th1}, we have
				\begin{equation}
			~\beta(\alpha,S_{\textbf{Q}})\geq \begin{cases}
						~	\alpha, \quad \quad&\alpha\in(0,\frac{1}{2}],   \\
				~	\frac{1}{2} , \quad \quad&\alpha\in(\frac{1}{2},\frac{3}{2}],  \\
				~	\frac{\alpha}{3}, \quad \quad&\alpha\in(\frac{3}{2},2],   \\
				~	\frac{2(\alpha-1)}{3} , \quad \quad&\alpha\in(2,4],  \\
				~	\alpha-2, \quad  \quad& \alpha\in(4,5].
			\end{cases}
		\end{equation}
	  \item[(2)] For $\textbf{Q}$ satisfying ${\rm(b)}$ in Theorem \ref{th1}, we have
	  		\begin{equation}
	  	~\beta(\alpha,S_{\textbf{Q}})\geq \begin{cases}
	  				~	\alpha, \quad \quad&\alpha\in(0,\frac{1}{2}],   \\
	  		~	\frac{1}{2} , \quad \quad&\alpha\in(\frac{1}{2},\frac{3}{2}],  \\
	  		~	\frac{\alpha}{3}, \quad \quad&\alpha\in(\frac{3}{2},\frac{9}{5}],   \\
	  		~	\frac{3(\alpha-1)}{4}, \quad  \quad& \alpha\in(\frac{9}{5},5].
	  	\end{cases}
	  \end{equation}
  	  \item[(3)] For $\textbf{Q}$ satisfying ${\rm(c)}$ in Theorem \ref{th1}, we have
  \begin{equation}
  	~\beta(\alpha,S_{\textbf{Q}})\geq \begin{cases}
  			 ~	\alpha, \quad \quad&\alpha\in(0,\frac{1}{2}],   \\
  		~	\frac{1}{2} , \quad \quad&\alpha\in(\frac{1}{2},\frac{3}{2}],  \\
  		~	\frac{\alpha}{3}, \quad \quad&\alpha\in(\frac{3}{2},3],   \\
  		~	\alpha-2, \quad  \quad& \alpha\in(3,5].
  	\end{cases}
  \end{equation}
  	  \item[(4)] For $\textbf{Q}$ satisfying ${\rm(d)}$ in Theorem \ref{th1}, we have
\begin{equation}
	~\beta(\alpha,S_{\textbf{Q}})\geq \begin{cases}
				~	\alpha, \quad \quad&\alpha\in(0,\frac{1}{2}],   \\
		~	\frac{1}{2} , \quad \quad&\alpha\in(\frac{1}{2},1],  \\
		~	\frac{\alpha}{2}, \quad \quad&\alpha\in(1,4],   \\
		~	\alpha-2, \quad  \quad& \alpha\in(4,5].
	\end{cases}
\end{equation}
\item[(5)] Suppose that $\textbf{Q}$ satisfies ${\rm(e)}$ in Theorem \ref{th1}. If $\textbf{Q}$ satisfies $\mathfrak{d}_{2,1}(\textbf{Q})=\mathfrak{d}_{1,2}(\textbf{Q})=1$, then
\begin{equation}\label{th2e5}
	~\beta(\alpha,S_{\textbf{Q}})\geq \begin{cases}
		~	\alpha, \quad \quad&\alpha\in(0,\frac{1}{2}],   \\
		~	\frac{1}{2} , \quad \quad&\alpha\in(\frac{1}{2},1],  \\
		~   \frac{\alpha}{2},  \quad\quad&\alpha\in(1,2],\\
		~	\frac{2\alpha-1}{3}, \quad  \quad& \alpha\in(2,5].
	\end{cases}
\end{equation}
If $\textbf{Q}$ satisfies $\mathfrak{d}_{2,1}(\textbf{Q})=1$ and  $\mathfrak{d}_{1,2}(\textbf{Q})=0$, then
\begin{equation}\label{th2e6}
	~\beta(\alpha,S_{\textbf{Q}})\geq \begin{cases}
		~	\alpha, \quad \quad&\alpha\in(0,\frac{1}{2}],   \\
		~	\frac{1}{2} , \quad \quad&\alpha\in(\frac{1}{2},1],  \\
		~   \frac{\alpha}{2},  \quad\quad&\alpha\in(1,4],\\
		~	\alpha-2, \quad  \quad& \alpha\in(4,5].
	\end{cases}
\end{equation}
If $\textbf{Q}$ satisfies $\mathfrak{d}_{2,1}(\textbf{Q})=0$ and  $\mathfrak{d}_{1,2}(\textbf{Q})=1$, then
\begin{equation}\label{th2e7}
	~\beta(\alpha,S_{\textbf{Q}})\geq \begin{cases}
		~	\alpha, \quad \quad&\alpha\in(0,\frac{1}{2}],   \\
		~	\frac{1}{2} , \quad \quad&\alpha\in(\frac{1}{2},1],  \\
		~   \frac{\alpha}{2},  \quad\quad&\alpha\in(1,3],\\
		~	\frac{3(\alpha-1)}{4}, \quad  \quad& \alpha\in(3,5].
	\end{cases}
\end{equation}
If $\textbf{Q}$ satisfies $\mathfrak{d}_{2,1}(\textbf{Q})=\mathfrak{d}_{1,2}(\textbf{Q})=0$, then
\begin{equation}\label{th2e8}
	~\beta(\alpha,S_{\textbf{Q}})\geq \begin{cases}
		~	\alpha, \quad \quad&\alpha\in(0,\frac{1}{2}],   \\
		~	\frac{1}{2} , \quad \quad&\alpha\in(\frac{1}{2},1],  \\
		~   \frac{\alpha}{2},  \quad\quad&\alpha\in(1,4],\\
		~	\alpha-2, \quad  \quad& \alpha\in(4,5].
	\end{cases}
\end{equation}
	\end{itemize}	
\end{theorem}

Note (\ref{W1}), we immediately have the following weighted restriction results.

\begin{theorem}\label{th3}
	Let $\alpha \in (0,5]$. Then
	\begin{itemize}
		\item[(1)] For $\textbf{Q}$ satisfying ${\rm(a)}$ in Theorem \ref{th1}, we have
	\begin{equation}
		~s(\alpha,S_{\textbf{Q}}) \leq \begin{cases}
			~	0, \quad \quad&\alpha\in(0,\frac{1}{2}],   \\
			~	\frac{2\alpha-1}{4}, \quad \quad&\alpha\in(\frac{1}{2},\frac{3}{2}],  \\
			~	\frac{\alpha}{3}, \quad \quad&\alpha\in(\frac{3}{2},2],   \\
~	\frac{\alpha+2}{6} , \quad \quad& \alpha\in(2,4],  \\
~	1, \quad  \quad& \alpha\in(4,5].
		\end{cases}
	\end{equation}
		\item[(2)] For $\textbf{Q}$ satisfying ${\rm(b)}$ in Theorem \ref{th1}, we have
		\begin{equation}
	~s(\alpha,S_{\textbf{Q}}) \leq \begin{cases}
					~	0, \quad \quad&\alpha\in(0,\frac{1}{2}],   \\
		~	\frac{2\alpha-1}{4}, \quad \quad&\alpha\in(\frac{1}{2},\frac{3}{2}],  \\
		~	\frac{\alpha}{3}, \quad \quad& \alpha\in(\frac{3}{2},\frac{9}{5}],   \\		~	\frac{\alpha+3}{8} , \quad \quad&  \alpha\in(\frac{9}{5},5].
	\end{cases}
\end{equation}
		\item[(3)] For $\textbf{Q}$ satisfying ${\rm(c)}$ in Theorem \ref{th1}, we have
\begin{equation}
	~s(\alpha,S_{\textbf{Q}}) \leq \begin{cases}
			~	0, \quad \quad&\alpha\in(0,\frac{1}{2}],   \\
		~	\frac{2\alpha-1}{4}, \quad \quad&\alpha\in(\frac{1}{2},\frac{3}{2}],  \\
		~	\frac{\alpha}{3}, \quad \quad& \alpha\in(\frac{3}{2},3],   \\
			~	1, \quad  \quad& \alpha\in(3,5].
	\end{cases}
\end{equation}
		\item[(4)] For $\textbf{Q}$ satisfying ${\rm(d)}$ in Theorem \ref{th1}, we have
		\begin{equation}
~s(\alpha,S_{\textbf{Q}}) \leq \begin{cases}
		~	0, \quad \quad&\alpha\in(0,\frac{1}{2}],   \\
	~	\frac{2\alpha-1}{4}, \quad \quad&\alpha\in(\frac{1}{2},1],  \\
	~	\frac{\alpha}{4}, \quad \quad& \alpha\in(1,4],   \\
	~	1, \quad  \quad&  \alpha\in(4,5].
\end{cases}
		\end{equation}
\item[(5)] Suppose that $\textbf{Q}$ satisfies ${\rm(e)}$ in Theorem \ref{th1}. If $\textbf{Q}$ satisfies  $\mathfrak{d}_{2,1}(\textbf{Q})=\mathfrak{d}_{1,2}(\textbf{Q})=1$,  then
\begin{equation}\label{th3e5}
~s(\alpha,S_{\textbf{Q}}) \leq \begin{cases}
			~	0, \quad \quad&\alpha\in(0,\frac{1}{2}],   \\
	~	\frac{2\alpha-1}{4}, \quad \quad&\alpha\in(\frac{1}{2},1],  \\
	~	\frac{\alpha}{4}, \quad \quad& \alpha\in(1,2],  \\
	~	\frac{\alpha+1}{6} , \quad \quad& \alpha\in(2,5].
\end{cases}
\end{equation}
If $\textbf{Q}$ satisfies  $\mathfrak{d}_{2,1}(\textbf{Q})=1$ and  $\mathfrak{d}_{1,2}(\textbf{Q})=0$, then
\begin{equation}\label{th3e6}
~s(\alpha,S_{\textbf{Q}}) \leq \begin{cases}
				~	0, \quad \quad&\alpha\in(0,\frac{1}{2}],   \\
	~	\frac{2\alpha-1}{4}, \quad \quad&\alpha\in(\frac{1}{2},1],  \\
	~	\frac{\alpha}{4}, \quad \quad& \alpha\in(1,4],  \\
	~	1, \quad  \quad& \alpha\in(4,5].
\end{cases}
\end{equation}
If $\textbf{Q}$ satisfies  $\mathfrak{d}_{2,1}(\textbf{Q})=0$ and  $\mathfrak{d}_{1,2}(\textbf{Q})=1$,  then
\begin{equation}\label{th3e7}
~s(\alpha,S_{\textbf{Q}}) \leq \begin{cases}
				~	0, \quad \quad&\alpha\in(0,\frac{1}{2}],   \\
	~	\frac{2\alpha-1}{4}, \quad \quad&\alpha\in(\frac{1}{2},1],  \\
	~	\frac{\alpha}{4}, \quad \quad& \alpha\in(1,3],  \\
	~	\frac{\alpha+3}{8} , \quad \quad&  \alpha\in(3,5].
\end{cases}
\end{equation}
If $\textbf{Q}$ satisfies  $\mathfrak{d}_{2,1}(\textbf{Q})=\mathfrak{d}_{1,2}(\textbf{Q})=0$, then
\begin{equation}\label{th3e8}
~s(\alpha,S_{\textbf{Q}}) \leq \begin{cases}
				~	0, \quad \quad&\alpha\in(0,\frac{1}{2}],   \\
	~	\frac{2\alpha-1}{4}, \quad \quad&\alpha\in(\frac{1}{2},1],  \\
	~	\frac{\alpha}{4}, \quad \quad& \alpha\in(1,4],   \\
	~	1, \quad  \quad&  \alpha\in(4,5].
\end{cases}
\end{equation}
	\end{itemize}
\end{theorem}

For the standard parabolic case ($n=1$ and $\textbf{Q}=\xi_1^2+...+\xi_d^2$), Du and Zhang \cite{DZ} gained the best bound in (\ref{qwe}) so far with $s=\frac{\alpha}{2(d+1)}$. Now we give a brief sketch of their proof. By locally constant property (see Section 2.2), (\ref{qwe}) can be reduced to
\begin{equation}\label{s1 eq1}
\|E^{\textbf{Q}}f\|_{L^2(X)} \lesssim R^{\frac{\alpha}{2(d+1)}} \|f\|_{L^2}, 
\end{equation}
where $X$ is an union of lattice unit cubes in $B^{d+1}(0,R)$ satisfying appropriate sparse property. They introduced a new parameter $\gamma$ to measure the sparse property of $X$ (see (\ref{sign000}) for the precise definition),  and established a more refined estimate:
\begin{equation}\label{s1 eq2}
	\|E^{\textbf{Q}}f\|_{L^2(X)} \lesssim \gamma^{\frac{1}{n+1}} R^{\frac{\alpha}{2(d+1)}} \|f\|_{L^2},  
\end{equation}
where they called it the fractal $L^2$ restriction estimate. In fact, they also established a more refined estimate than (\ref{s1 eq2}), see \cite[Theorem 1.6]{DZ}. The estimate (\ref{s1 eq2}) is very strong, and it implies not only (\ref{s1 eq1}), but also the almost sharp pointwise convergence of the solution to the free Schr\"odinger equation. The main method to prove (\ref{s1 eq2}) is the broad-narrow analysis developed by Bourgain and Guth \cite{BG11}. Du and Zhang separated the original problem into two cases with the  dimension $d$ as a division: the essentially $(d+1)$-dimensional case (the broad case), and the $d$-dimensional case (the narrow case). For the broad case, the multilinear restriction estimate can offer a good bound in (\ref{s1 eq2}). For the narrow case, the interaction between square root cancellation produced by the lower dimensional $\ell^2$ decoupling and the  additional parameter $\gamma$ makes the scale induction close. If we separate this  problem with general dimension $k$ as a division ($1\leq k \leq d$) and carry out their argument, we find that  the corresponding fractal $L^2$ restriction estimate is best only when $k=d$, which depends on the geometry of the paraboloid. For the conical case, Harris \cite{H19} separated this problem with the dimension $d-1$ as a division. For the hyperbolic case,  Barron, Erdo\v{g}an and Harris \cite{BEH20} separated this problem with general dimension $k$ as a division which relies on the different $\ell^2$ decoupling inequalities for the hyperboloids with different signatures (see Theorem \ref{dec 2}). 

For the cases of high co-dimensions, things turn out more complicated. Firstly, there may be many different kinds of quadratic forms as the increase of co-dimension $n$.
Let us take $d=3$ as an example. When $n=1$, there exist only two kinds of quadratic forms: $\textbf{Q}=\xi_1^2+\xi_2^2+\xi_3^2$ corresponding to the paraboloid, and  $\textbf{Q}=\xi_1^2+\xi_2^2-\xi_3^2$ corresponding to the hyperboloid. When $n=2$, we can see that there are many kinds of quadratic forms in Theorem \ref{th1}. Next, degenerate quadratic forms possess a few strange characteristics. Here we roughly call that $\textbf{Q}$ is degenerate if the values of   $\mathfrak{d}_{d',n'}(\textbf{Q})$ are small (see \cite{GOZZK} for the precise definition), since the $\mathfrak{d}_{d',n'}(\textbf{Q})$ reflects the ``curvature" of the associated quadratic surface $S_{\textbf{Q}}$.  For instance, the lower dimensional decoupling is not always better than the original dimensional decoupling for such quadratic forms, and we will observe this phenomenon in Corollary \ref{dqf3}. On the other hand, such quadratic forms $\textbf{Q}$ possess some special scaling structures. We now explain this point by discussing the relevant restriction problems. The goal of the restriction problems on $\textbf{Q}$ is to find the optimal ranges of $p$ and $q$ such that the inequality 
\begin{equation}\label{rs}
	\|E^{\textbf{Q}}f\|_{L^q(\mathbb{R}^{d+n})} \lesssim  \|f\|_{L^p(\mathbb{R}^d)}
\end{equation}
holds for every function $f \in L^p(\mathbb{R}^d)$. When $n=1$, we can use the classical  Knapp example $\chi_{[0,R^{-1}]^d}$ to attain the necessary condition of (\ref{rs}). When $n \geq 2$, the same example may not to lead to the optimal necessary condition. Readers can see Table 1 for the quadratic forms $\textbf{Q}$ in Theorem \ref{th1}. It is interesting to observe that the quadratic surfaces $S_\textbf{Q}$ remain invariant if we scale the associated examples in Table 1 to $\chi_{[0,1]^3}$. Besides this, we also point out that two quadratic forms $\textbf{Q}$ and $\textbf{Q}'$ may have different restriction results even if $ \mathfrak{d}_{d',n'}(\textbf{Q})= \mathfrak{d}_{d',n'}(\textbf{Q}')$ for all  $0 \leq d' \leq d$ and $0\leq n'\leq n$, such as  $(\xi_1\xi_2,\xi_2^2+\xi_1\xi_3)$ and $ (\xi_1\xi_2,\xi_1^2\pm\xi^2_3)$. In order to utilize every quadratic form $\textbf{Q}$'s scaling structure and degenerate characteristic more sufficiently, we will use the broad-narrow analysis with different dimensions as divisions, together with adaptive $\ell^2$ decoupling, to attain a few weighted restriction estimates in Theorem \ref{th3}.

\begin{table}[htbp]
	
	\begin{tabular}{|c|c|c|}
		\hline
		\rule{0pt}{15pt}
		$\textbf{Q}$ & \makecell[c]{  Examples: characteristic functions  \\ with support   }& \makecell[c]{Necessary conditions \\ of (\ref{rs})}  \\
		\hline
		\rule{0pt}{15pt}
		$(\xi_1^2,\xi_2^2+\xi_1\xi_3)$  & $[0,R^{-1}] \times [0,R^{-1/2}]\times [0,1]$ & $q > 4,~q\geq 3p'$ \\
		\hline \rule{0pt}{15pt}
		$(\xi_1^2,\xi_2\xi_3)$ & $[0,R^{-1}] \times [0,1]\times [0,1]$&  $q > 4,~q\geq 3p'$\\
		\hline \rule{0pt}{15pt}
		$(\xi_1^2,\xi_2^2+\xi_3^2)$ & $[0,R^{-1}] \times [0,1]\times [0,1]$ &  $q > 4,~q\geq 3p'$  \\
		\hline \rule{0pt}{15pt}
		$ (\xi_1 \xi_2,\xi_1\xi_3)$ & $[0,R^{-1}] \times [0,1]\times [0,1]$& $q > 4,~q\geq 3p'$  \\
		\hline \rule{0pt}{15pt}
		$(\xi_1\xi_2,\xi_2^2+\xi_1\xi_3)$ & $[0,R^{-1}] \times [0,R^{-1/2}]\times [0,1]$ &  $q > 11/3,~q\geq 8p'/3$  \\
		\hline \rule{0pt}{15pt}
		$ (\xi_1\xi_2,\xi_1^2\pm\xi^2_3)$ & $[0,R^{-1}] \times [0,1]\times [0,R^{-1}]$&  $q > 7/2,~q\geq 5p'/2$ \\
		\hline \rule{0pt}{15pt}
		satisfies the (CM) condition & $[0,R^{-1}] \times [0,R^{-1}] \times [0,R^{-1}] $&  $q > 10/3,~q\geq 7p'/3$ \\
		\hline
	\end{tabular}
	\vskip0.2cm
	\centering
	\caption{The necessary conditions and corresponding examples for the quadratic forms in Theorem \ref{th1}. See \cite[Section 3]{GO}. }
\end{table}

However, this is not sufficient to prove all bounds in Theorem \ref{th3}. To gain better  bounds for small $\alpha$, we will build the almost sharp Stein-Tomas type's inequalities for the quadratic forms in Theorem \ref{th1}. The Stein-Tomas inequalities for hypersurfaces of nonzero Gaussian curvature were built by the interpolation of analytic families of operators. Then Christ \cite{C82,C85} and Mockenhaupt \cite{M96} extended this argument to the case of non-degenerate surfaces of high co-dimensions. In particular, they attained the sharp Stein-Tomas type's inequalities for $\textbf{Q}$ satisfying ${\rm(e)}$ in Theorem \ref{th1}. Nevertheless, the same proof is invalid for degenerate quadratic forms. As a substitute, we will use the broad-narrow analysis, combined with the bilinear restriction estimates for quadratic forms built in this paper, to get the almost sharp Stein-Tomas type's inequalities for such  quadratic forms.

\vskip0.3cm

\noindent \textbf{Outline of paper.} In Section 2, we present some important tools, which will be used throughout the rest of the paper. Our argument relies on the Du-Zhang method of \cite{DZ} from Section 3 to Section 6. In Section 3, we use the multilinear restriction estimates  and lower dimensional $\ell^2$ decoupling for quadratic forms to obtain Proposition \ref{sec2 pro2}, from which we attain a few weighted restriction results in Corollary \ref{cor1}. In Section 4, we use the $\ell^2$ decoupling for quadratic forms directly to obtain Proposition \ref{sec3 pro1}, from which we attain a few weighted restriction results in Corollary \ref{cor2}. In Section 5, we establish one class of $k$-linear restriction estimates and Kakeya estimates for quadratic forms by the Brascamp-Lieb inequality. In Section 6, we use the bilinear restriction estimates built in the previous section and the adaptive $\ell^2$ decoupling to obtain a few weighted restriction results in Proposition \ref{cor3}. The almost sharp Stein-Tomas type's inequalities for the quadratic forms in Theorem \ref{th1} are built in Section 7. As a direct result, we obtain a few weighted restriction results in Corollary \ref{cor4}. Combining Corollary \ref{cor1}, Corollary \ref{cor2}, Proposition \ref{cor3} and Corollary \ref{cor4}, we have showed Theorem \ref{th3}. Finally, in Appendix A, we build the refined Strichartz estimates for quadratic forms, which can provide another alternative proof of Corollary \ref{buchong2}. And in Appendix B, we show Theorem \ref{th1} and Theorem 1.2 from \cite{GO} are equivalent.

\vskip0.3cm

\noindent \textbf{Notation.} If $X$ is a finite set, we use $\# X$ to denote its cardinality. If $X$ is a measurable set, we use $|X|$ to denote its Lebesgue measure. We use $B^m(c,r)$ to represent a ball centered at $c$ with radius $r$ in $\mathbb{R}^m$. We abbreviate $B^m(c,r)$ to $B(c,r)$ if $\mathbb{R}^m$ is clear in the context. If $p \geq 1$, we use $p'$ to  denote the dual exponent of $p$. We write $A\lesssim_\epsilon B$ to
mean that there exists a constant $C$ depending on $\epsilon$ such that $A\leq CB$. Moreover, $A \sim B$ means $A \lesssim B$ and $A\gtrsim B$. We use $\xi^\top$ to denote the transpose of $\xi$.  Define $e(b):=e^{2\pi i b}$ for each $b \in \mathbb{R}$.

Let $\delta\in(0,1)$ be a dyadic number. We denote by $\mathcal{P}(Q,\delta)$ the dyadic cubes of side length $\delta$ in $Q$ for every dyadic cube $Q \in [0,1]^d$. Let $\mathcal{P}(\delta)$ be the partition of $[0,1]^d$ into the dyadic cubes of side length $\delta$. If $W$ is a cube with side length $l(W)$, we use $c\cdot W$ to denote the cube of side length $c \cdot l(W)$ and of the same center as $W$. If $S$ is a surface, we use $N_\delta(S)$ to denote $\delta$-neighborhood of $S$. 

For a ball $B=B(c_B, r_B)$ with center $c_B$ and radius $r_B$, define an associated weight
\begin{equation}\label{sig1}
	w_B(\cdot):=\Big(   1+ \frac{|\cdot-c_B|}{r_B} \Big)^{-N}, 
\end{equation}
with a large constant $N>0$. Then we define averaged integrals:
$$  \|f\|_{L^p_{\text{avg}}(B)}:=\Big(   \frac{1}{|B|} \int_B |f(x)|^pdx  \Big)^{\frac{1}{p}}~\text{~~and~~}~ \|f\|_{L^p_{\text{avg}}(w_B)}:=\Big(   \frac{1}{|B|} \int |f(x)|^p w_B(x)dx \Big)^{\frac{1}{p}}.   $$

\section{Preliminaries}

We summarize some important results in this section.

\vskip0.3cm

\noindent 2.1. \textbf{Basic positive results.}
\begin{lemma}[\cite{S97,E04}]\label{Lem2.1}
	Let $m\geq 1$ and $\alpha \in (0,m)$. Let $\sigma$ be a finite measure on $\mathbb{R}^m$ with compact support such that
	$$|\widehat{\sigma}(\xi)|\lesssim|\xi|^{-a}~~~\quad ~\text{and}~\quad ~~~\sigma(B(x,r))\lesssim r^b,$$
	for $a,b\in(0,m)$. Then for any $\mu\in\mathcal{M}(\mathbb{R}^m)$ with $I_{\alpha}(\mu)<\infty$, we have
	$$\int|\widehat{\mu}(R\xi)|^2d\sigma(\xi)\lesssim R^{-\max\{\min\{a,\alpha\},\alpha-m+b\}}.$$
\end{lemma}
In this setting, for the quadratic forms $\textbf{Q}$ satisfying the (CM) condition, it was showed in \cite[Section 4]{GO} that
\begin{equation}\label{add2}
|E^{\textbf{Q}}1(x,y)|\lesssim(1+|(x,y)|)^{-1/2}.
\end{equation}
In fact, for every quadratic form listed in Theorem \ref{th1},  (\ref{add2}) always holds. Suppose that the quadratic form $\textbf{Q}=(Q_1,Q_2)$ satisfies the condition of Theorem \ref{th1}. Write
\begin{align*}
	 E^{\textbf{Q}}1(x,y)&=\int_{[0,1]^3} e\Big(x \cdot \xi^\top + y\cdot (Q_1,Q_2)^\top\Big) d\xi  \\
	 &  =\int_{[0,1]^3} e\Big(x \cdot \xi^\top + y M_1 \cdot M_1^{-1}(Q_1,Q_2)^\top\Big) d\xi.
\end{align*}
Here $M_1$ denotes the rotation matrix, which is defined by
\begin{equation}
	M_1:= \begin{pmatrix}
		\cos \theta & \sin\theta  \\
		-\sin \theta & \cos \theta 
	\end{pmatrix},
\end{equation}
such that $y M_1=(y_1,0)$. Therefore, one obtains 
\begin{align*}
E^{\textbf{Q}}1(x,y) &=\int_{[0,1]^3} e\Big[x \cdot \xi^\top + y_1  \big(Q_1(\xi)\cos \theta +Q_2(\xi)\sin\theta \big)\Big] d\xi \\ 
&= \int_{[0,1]^3} e\Big[ |(x,y_1)|  \Big( \frac{x}{|(x,y_1)|} \cdot \xi^\top + \frac{y_1}{|(x,y_1)|}  \big(Q_1(\xi)\cos \theta +Q_2(\xi)\sin\theta \big)\Big)\Big] d\xi .
\end{align*}
Define 
$$\tilde{Q}(\xi)=Q_1(\xi)\cos \theta +Q_2(\xi)\sin\theta,$$ 
and
$$  \phi(\xi)=  \frac{x}{|(x,y_1)|} \cdot \xi^\top + \frac{y_1}{|(x,y_1)|}  \tilde{Q}(\xi), \quad \quad (x,y_1)\neq 0. $$
Note that 
$$  \bigtriangledown_\xi \phi(\xi)= \frac{x}{|(x,y_1)|}+ \frac{y_1}{|(x,y_1)|}  \bigtriangledown_\xi \tilde{Q}(\xi),$$
and 
$$ \bigtriangledown^2_\xi \phi(\xi)=  \frac{y_1}{|(x,y_1)|}  \bigtriangledown^2_\xi \tilde{Q}(\xi).   $$
Now we fix $(x,y_1)\neq 0$. If $|y_1|/|(x,y_1)| \ll 1$, then $|x|/|(x,y_1)| \sim 1$. On the other hand, we know that $|\bigtriangledown_\xi \tilde{Q}(\xi)|$ is bounded in $[0,1]^3$. Thus we get $| \bigtriangledown_\xi \phi(\xi)|\gtrsim 1$. Using van der Corput's lemma, such as Theorem 7.2 in \cite{CCW99}, we have
\begin{equation*}
|E^{\textbf{Q}}1(x,y)|\lesssim(1+|(x,y)|)^{-1}.
\end{equation*}
If $|y_1|/|(x,y_1)| \sim 1$, then $|\bigtriangledown^2_\xi \phi(\xi)| \sim    |\bigtriangledown^2_\xi \tilde{Q}(\xi)|. $ Firstly, we assume that $\textbf{Q}= (\xi_1^2,\xi_2^2+\xi_1\xi_3)$, and then
\begin{equation}
	\bigtriangledown^2_\xi \tilde{Q}(\xi)= \begin{pmatrix}
		2\cos \theta & 0 &\sin\theta  \\
	   0 &	2\sin \theta & 0  \\
	   \sin \theta & 0 & 0 
	\end{pmatrix}.
\end{equation}
This implies that $|\partial^\beta \tilde{Q}| \gtrsim 1$ for $\beta=(1,1)$ or $\beta=(2,2)$. Then using Theorem 7.2 in \cite{CCW99} again, we get  
\begin{equation*}
	|E^{\textbf{Q}}1(x,y)|\lesssim(1+|(x,y)|)^{-1/2}.
\end{equation*}
So we have proved (\ref{add2}) for $\textbf{Q}= (\xi_1^2,\xi_2^2+\xi_1\xi_3)$. On the cases (a)-(d) in Theorem \ref{th1}, we can obtain (\ref{add2}) through the same argument. Let $\sigma$ be the surface measure on $S_{\textbf{Q}}$, by Lemma \ref{Lem2.1}, then one concludes
\begin{align*}
	\begin{split}
		\beta(\alpha,S_{\textbf{Q}})\geq\left\{
		\begin{array}{ll}
			\alpha,                    \quad\quad   &\alpha\in(0,\frac{1}{2}],\\
			\frac{1}{2},             \quad\quad              &\alpha\in(\frac{1}{2},\frac{5}{2}], \\
			\alpha-2,              \quad\quad        &\alpha\in(\frac{5}{2},5].
		\end{array}
		\right.
	\end{split}
\end{align*}
On the other hand, by (\ref{W1}), we have
\begin{align*}
	\begin{split}
		s(\alpha,S_{\textbf{Q}})\leq\left\{
		\begin{array}{ll}
			0,                    \quad\quad   &\alpha\in(0,\frac{1}{2}],\\
			\frac{2\alpha-1}{4},             \quad\quad              &\alpha\in(\frac{1}{2},\frac{5}{2}], \\
		1,              \quad\quad        &\alpha\in(\frac{5}{2},5].
		\end{array}
		\right.
	\end{split}
\end{align*}
These two results give the bounds in Theorem \ref{th2} and Theorem \ref{th3} for small $\alpha$, respectively.

\vskip0.3cm

\noindent 2.2. \textbf{Locally constant property.} 
\begin{lemma}[\cite{GWZ}, Lemma 6.1, Lemma 6.2]\label{lcp}
Let $m \geq 1$, $\theta$ be a compact symmetric convex set centered at $c(\theta) \in \mathbb{R}^m$. Let $T_\theta:=\{ x:|x\cdot (y-c(\theta))| \leq 1,\forall y \in \theta \}$, and $T$ be a translated copy of $T_\theta$. If a function $f$ is Fourier supported in $\theta$, then there exists a positive function $\eta_{T_\theta}$ satisfying the following properties:
\begin{itemize}
	\item[(1)] $\eta_{T_\theta}$ is essentially supported on $10T_\theta$ and rapidly decaying away from it, and $\|\eta_{T_\theta}\|_{L^1} \lesssim 1$. Here rapidly decaying means that for any integer $N>0$, there exists a constant $C_N$ such that $\eta_{T_\theta}(x)\leq C_N (1+n(x,10T_\theta))^{-N}$, where $n(x,10T_\theta)$ is the smallest positive integer $n$ such that $x \in 10nT_\theta$. 
	\item[(2)] Let $c_T:= \max_{x \in T} |f(x)|$, then we have
	\begin{equation}
     |f|\leq \sum_{T\pxx T_\theta} c_T \chi_T \leq |f| \ast \eta_{T_\theta},
	\end{equation}
where the sum $\sum_{T\pxx T_\theta} $ is over a finitely overlapping cover $\{T\}$ of $\mathbb{R}^m$ with $T\pxx T_\theta$.
	\item[(3)] For any $1 \leq p< \infty $, one has
	\begin{equation}
		\int_T  (|f| \ast \eta_{T_\theta})^p \lesssim_p \int |f|^p w_T,
	\end{equation}
where $w_T=1$ on $10T$, and rapidly decaying away from it.
\end{itemize}
\end{lemma}

Locally constant property says that if a function $f$ has compact Fourier support $\theta$, then we can view $|f|$ as constant on every dual set $T$ in the sense of integral average. More precisely, for any $1 \leq p,q< \infty $,  we have
	\begin{equation}
	\|f\|_{L^p_{{\rm avg}}(w_T)} \sim 	\|f\|_{L^q_{{\rm avg}}(w_T)}.
\end{equation}
This does not imply that $|f|$ is essentially constant on $T$, i.e. $|f| \sim c$ for some constant $c>0$. A simple example is the function $f$ such that $\widehat{f}=\chi_\theta$.  Nevertheless, if we take
$$   \tilde{w}_T(x):=  |T|^{-1}  (1+n(x,10T_\theta))^{-100m},  $$
which is essentially constant on every $T$, then by
$$   |f|\leq |f| \ast \eta_{T_\theta} \lesssim|f| \ast \tilde{w}_T,    $$
one gets $|f| \ast \tilde{w}_T$  is essentially constant on every $T$.
Hence we can deal with $|f|$ as constant on every $T$ in the sense of convolution.

\vskip0.3cm

\noindent 2.3. \textbf{Wave packet decomposition.} Let $d \geq 1$ and $n \geq 1$. Let $\textbf{Q}=(Q_1,...,Q_n)$ be a collection of quadratic forms in $d$ variables. Fix a scale $R >1$. Cover $B^d(0,1)$ by finitely overlapping balls $\theta$ with radius $R^{-1/2}$ and center $c(\theta)$. Cover $B^d(0,R)$ by finitely overlapping balls $\nu$ with radius $R^{1/2}$ and center $c(\nu)$. Using partition of unity, we have 
$$f=\sum_{\theta,\nu} f_{\theta,\nu},$$
where each $f_{\theta,\nu}$ is of Fourier supported in $\theta$, and essentially supported in $\nu$. The functions $f_{\theta,\nu}$ are approximately orthogonal:
$$    \Big\|   \sum_{\theta,\nu} f_{\theta,\nu}\Big\|_{L^2}^2 \sim  \sum_{\theta,\nu} \left\|    f_{\theta,\nu}\right\|_{L^2}^2.  $$
On the other hand, we have
$$E^{\textbf{Q}}f= \sum_{\theta,\nu} E^{\textbf{Q}}f_{\theta,\nu},$$
and each $E^{\textbf{Q}}f_{\theta,\nu}$ restricted on $B^{d+n}_R$ is essentially supported on a rectangular box $\Box_{\theta,\nu}$:
$$\Box_{\theta,\nu}:=\Big\{    (x,y) \in B^{d+n}_R  : \Big|x+\sum_{j=1}^n y_j \nabla Q_j(c(\theta))-c(\nu)    \Big|\leq R^{\frac{1}{2}+\delta}  \Big\},$$
where $\delta$ is a small parameter. For more details on this wave packet decomposition see \cite{BLL17}.

\vskip0.3cm

\noindent 2.4. \textbf{Decoupling.} We introduce the decoupling theorem for quadratic forms in \cite{GOZZK}. 
 Let $q,p \geq 2$, and $0< \delta<1$ be a dyadic number. Given a function $f$ supported on $[0,1]^d$, let $f=\sum_\tau f_\tau$, where each $f_\tau$ is supported in a dyadic cube $\tau$ of side length $\delta$. Define $\Gamma_{q,p}^d(\textbf{Q})$ be the smallest constant $\Gamma$ such that
\begin{equation}
	\|E^{\textbf{Q}}f\|_{L^p(\mathbb{R}^{d+n})} \lesssim \delta^{-\Gamma-\epsilon} \Big(   \sum_\tau \|E^{\textbf{Q}}f_\tau\|_{L^p(\mathbb{R}^{d+n})}^q  \Big)^{\frac{1}{q}}, \quad \forall~\epsilon>0, ~\delta<1.
\end{equation}

\begin{theorem}[\cite{GOZZK}, Theorem 1.1]\label{dqf} 
Let $d \geq 1$ and $n \geq 1$. Let $\textbf{Q}=(Q_1,...,Q_n)$ be a collection of quadratic forms in $d$ variables. For $2 \leq q \leq p <\infty$, we have
\begin{equation}
\Gamma_{q,p}^d(\textbf{Q})= \max_{0\leq d'\leq d} \max_{0\leq n' \leq n} \Big\{  d'\Big(1-\frac{1}{p}-\frac{1}{q}\Big)- \mathfrak{d}_{d',n'}(\textbf{Q})\Big(\frac{1}{2}-\frac{1}{p}\Big) -\frac{2(n-n')}{p}    \Big\}.
\end{equation}
For $2 \leq p<q \leq \infty$, we have
\begin{equation}
	\Gamma_{q,p}^d(\textbf{Q})=\Gamma_{p,p}^d(\textbf{Q})+d\Big(\frac{1}{p}-\frac{1}{q}\Big).
\end{equation}
\end{theorem}

We say that $\textbf{Q}$ is strongly non-degenerate if
$$   \mathfrak{d}_{d',n'}(\textbf{Q}) \geq d'-\Big(1-\frac{n'}{n}\Big)d,   $$
for every $d'$ with $0\leq d'\leq d$ and every $n'$ with $0\leq n' \leq n$. We have the nice decoupling inequality for such $\textbf{Q}$, which can be viewed as a generalization for the case of the paraboloid. 

\begin{corollary}[\cite{GOZZK}, Corollary 1.3]\label{dqf2} 
$\textbf{Q}$ is strongly non-degenerate if and only if
\begin{equation}
	\Gamma_{2,p}^d(\textbf{Q})=\max\Big\{   0,d\Big(\frac{1}{2}-\frac{1}{p}\Big)-\frac{2n}{p}  \Big\}, \quad\quad 2\leq p<\infty. 
\end{equation}
\end{corollary}

Since we will only use the case $q=2$ in Theorem \ref{dqf} in the rest of the paper, we abbreviate  $\Gamma_p^d(\textbf{Q}):=\Gamma_{2,p}^d(\textbf{Q})$ for convenience.  
Now, combining with Theorem \ref{th1} and Theorem \ref{dqf}, we obtain the following corollary. 
\begin{corollary}[$\ell^2$ decoupling for surfaces of co-dimension 2 in $\mathbb{R}^5$]\label{dqf3}  Let $d=3, n=2$. Then
	
	\begin{itemize}
		\item[(1)] For $\textbf{Q}$ satisfying ${\rm(a)}$ and ${\rm(d)}$ in Theorem \ref{th1}, we have the sharp  decoupling exponent
		$$\Gamma_p^3(\textbf{Q})=\max\Big\{ \frac{1}{2}-\frac{1}{p},\frac{3}{2}-\frac{7}{p}   \Big\},$$
		with $\Gamma=1/3$ at the critical point $p=6$. We also have the sharp lower dimensional  decoupling exponent
		$$\Gamma_p^2(\textbf{Q})=\max\Big\{ \frac{1}{2}-\frac{1}{p},1-\frac{4}{p}   \Big\},$$
		with $\Gamma=1/3$ at the critical point $p=6$. 
		\item[(2)] For $\textbf{Q}$ satisfying ${\rm(b)}$ in Theorem \ref{th1}, we have the sharp decoupling exponent
		$$\Gamma_p^3(\textbf{Q})=\max\Big\{0, 1-\frac{4}{p},\frac{3}{2}-\frac{7}{p}   \Big\},$$
		with $\Gamma=0$ at the critical point $p=4$ and $\Gamma=1/3$ at the critical point $p=6$. We also have the sharp lower dimensional decoupling exponent
		$$\Gamma_p^2(\textbf{Q})=\max\Big\{ 0,1-\frac{4}{p}   \Big\},$$
	with $\Gamma=0$ at the critical point $p=4$. 
		\item[(3)] For $\textbf{Q}$ satisfying ${\rm(c)}$ in Theorem \ref{th1}, we have the sharp decoupling exponent 
		$$\Gamma_p^3(\textbf{Q})=\max\Big\{ 1-\frac{2}{p},\frac{3}{2}-\frac{7}{p}   \Big\},$$
	with $\Gamma=4/5$ at the critical point $p=10$. We also have the  sharp lower dimensional  decoupling exponent 
		$$\Gamma_p^2(\textbf{Q})=1-\frac{2}{p}.$$  
		\item[(4)] Suppose that $\textbf{Q}$ satisfies ${\rm(e)}$ in Theorem \ref{th1}. If $\textbf{Q}$ satisfies  $\mathfrak{d}_{2,1}(\textbf{Q})=\mathfrak{d}_{1,2}(\textbf{Q})=1$ which is strongly non-degenerate, we have the sharp decoupling exponent
		$$\Gamma_p^3(\textbf{Q})=\max\Big\{ 0,\frac{3}{2}-\frac{7}{p}   \Big\},$$
	with $\Gamma=0$ at the critical point $p=14/3$. We also have the  sharp  lower dimensional  decoupling exponent 
		$$\Gamma_p^2(\textbf{Q})=\max\Big\{ 0,1-\frac{6}{p}   \Big\},$$      
	with $\Gamma=0$ at the critical point $p=6$. 
	
	\noindent If $\textbf{Q}$ satisfies  $\mathfrak{d}_{2,1}(\textbf{Q})=1$ and  $\mathfrak{d}_{1,2}(\textbf{Q})=0$, we have the sharp decoupling exponent
		$$\Gamma_p^3(\textbf{Q})=\max\Big\{ \frac{1}{2}-\frac{1}{p},\frac{3}{2}-\frac{7}{p}   \Big\},$$
	with $\Gamma=1/3$ at the critical point $p=6$. We also have the  sharp  lower dimensional  decoupling exponent 
		$$\Gamma_p^2(\textbf{Q})=\max\Big\{ \frac{1}{2}-\frac{1}{p},1-\frac{6}{p}   \Big\},$$      
	with $\Gamma=2/5$ at the critical point $p=10$. 
	
	\noindent 	If $\textbf{Q}$ satisfies  $\mathfrak{d}_{2,1}(\textbf{Q})=0$ and  $\mathfrak{d}_{1,2}(\textbf{Q})=1$, we have the same decoupling exponents as in ${\rm(2)}$. 
	
	\noindent 	If $\textbf{Q}$ satisfies  $\mathfrak{d}_{2,1}(\textbf{Q})=\mathfrak{d}_{1,2}(\textbf{Q})=0$, we have the same decoupling exponents as in ${\rm(1)}$.
	\end{itemize}
\end{corollary}

Besides the $\ell^2$ decoupling in the above corollary, we will also use the $\ell^2$ decoupling for the paraboloid, hyperboloid, and flat manifold in Section 6. We list these results as below.

Let $M$ be a $(m-1) \times (m-1)$ diagonal matrix with all nonzero entries, and define
$$   \mathbb{H}:=\{   (\xi,\langle M\xi,\xi\rangle):\xi \in B^{m-1}(0,1)   \}.   $$
We use $s$ to denote the minimum of the number of positive and negative entries of $\mathbb{H}$. If $s=0$, we write $\mathbb{P}^{m-1}=\mathbb{H}$ corresponding to the paraboloid. If $s>0$, we write $\mathbb{H}_s^{m-1}=\mathbb{H}$ corresponding to the hyperboloid.

\begin{theorem} [$\ell^2$ decoupling for the paraboloid, \cite{BD15}]\label{dec 1}  
Let $F:\mathbb{R}^m \rightarrow \mathbb{C}$ be such that ${\rm supp}\widehat{F}\subset N_{R^{-1}}(\mathbb{P}^{m-1})$. Divide this neighborhood into slabs $\theta$ with $m-1$ long directions of length $R^{-1/2}$ and one short direction of length $R^{-1}$. Write $F=\sum_\theta F_\theta$, where $\widehat{F}_\theta=\widehat{F} \chi_\theta$. Then 
\begin{equation}
	\|F\|_{L^p(\mathbb{R}^m)} \lesssim R^{\beta(p)+\epsilon} \Big( \sum_{\theta} \|F_\theta\|_{L^p(\mathbb{R}^m)}^2  \Big)^{\frac{1}{2}},
\end{equation}
where $\beta(p)=0$ when $2\leq p \leq \frac{2(m+1)}{m-1}$ and $\beta(p)=\frac{m-1}{4}-\frac{m+1}{2p}$ when $\frac{2(m+1)}{m-1}\leq p \leq \infty$.
\end{theorem}

\begin{theorem} [$\ell^2$ decoupling for the hyperboloid, \cite{BD17b}]\label{dec 2}  
Let $F:\mathbb{R}^m \rightarrow \mathbb{C}$ be such that ${\rm supp}\widehat{F}\subset N_{R^{-1}}(\mathbb{H}_s^{m-1})$. Divide this neighborhood into slabs $\theta$ with $m-1$ long directions of length $R^{-1/2}$ and one short direction of length $R^{-1}$. Write $F=\sum_\theta F_\theta$, where $\widehat{F}_\theta=\widehat{F} \chi_\theta$. Then 
\begin{equation}
	\|F\|_{L^p(\mathbb{R}^m)} \lesssim R^{\beta(p)+\epsilon} \Big( \sum_{\theta} \|F_\theta\|_{L^p(\mathbb{R}^m)}^2  \Big)^{\frac{1}{2}},
\end{equation}
where $\beta(p)=(\frac{1}{4}-\frac{1}{2p})s$ when $2\leq p \leq \frac{2(m+1-s)}{m-1-s}$ and $\beta(p)=\frac{m-1}{4}-\frac{m+1}{2p}$ when $\frac{2(m+1-s)}{m-1-s}\leq p \leq \infty$.
\end{theorem}

Interpolating two trivial estimates for $p=2$ and $p=\infty$, we easily obtain the following flat decoupling.

\begin{theorem} [Flat decoupling]\label{dec 3}  
	Let $R$ be a rectangular box in $\mathbb{R}^m$, and $R_1,...,R_L$ be a partition of $R$ into congruent boxes that are translates of each other. Write $F=\sum_j F_{j}$, where $\widehat{F}_{j}=\widehat{F} \chi_{B_j}$. Then for any $2\leq p,q\leq \infty$, we have
	\begin{equation}
		\|F\|_{L^p(\mathbb{R}^m)} \lesssim L^{1-\frac{1}{p}-\frac{1}{q}}\Big( \sum_{j} \|F_j\|_{L^p(\mathbb{R}^m)}^q  \Big)^{\frac{1}{q}}.
	\end{equation}
\end{theorem}

\section{The Du-Zhang method for general quadratic forms}

From now on, we begin to prove Theorem \ref{th3}. In this section, by adapting the Du-Zhang method in \cite{DZ}, we will build the fractal $L^2$ restriction estimates for quadratic surfaces.  Our argument is essentially the same as in \cite{DZ} except two points. Note that the multilinear restriction estimates are strongly different between the paraboloid and general quadratic surfaces, we will adopt the algorithm in \cite[Proposition 5.6]{GOZZK} instead of the dichotomy in the argument of \cite[Section 3]{DZ}. As this strategy changes, in the narrow case, we need to use the lower dimensional decoupling when frequency support is clustered near sub-varieties (Lemma \ref{dec low var}). Unfortunately, this lemma can produce several new scales. To overcome this problem, we will adopt the matched inductive pattern for each different scale.

\begin{definition}
Let $\theta \in (0,1]$, and $K \in 2^{\mathbb{N}}$ be a dyadic integer. We say a subset $\mathcal{W} \subset \mathcal{P}(1/K)$ is $\theta$-uniform if for every non-zero polynomial $P$ in $d$ variables with real coefficients of degree $\leq d$, one has
$$    \# \{    W \in \mathcal{W}:2W \cap Z_P \neq \varnothing  \} \leq \theta \cdot \#\mathcal{W}.  $$
Here $Z_P$ denotes the zero set of $P$. 
\end{definition}

\begin{lemma}[Multilinear restriction estimates for quadratic forms, \cite{GOZZK}]\label{mul th} 
Let $K \in 2^{\mathbb{N}}$ be a dyadic integer, and $R$ be a scale with $ R\geq K$. Let $\theta >0$, and $\{W_j \}_{j=1}^J \subset \mathcal{P}(1/K)$ be a $\theta$-uniform set of cubes. Let $B \subset \mathbb{R}^{d+n}$ be a ball of radius $R^2$. Then for each $2 \leq p<\infty$, any $\epsilon>0$, and every $f_j$ with ${\rm supp}f_j \subset W_j$, we have
\begin{equation}
	\left\|   \prod_{j=1}^{J}   |E^{\textbf{Q}}f_j|^{\frac{1}{J}} \right\|_{L^p(B)} \lesssim_{\theta,K,\epsilon} R^{\gamma(p,\theta,\textbf{Q})+\epsilon} \prod_{j=1}^{J}\|f_j\|^{\frac{1}{J}}_{L^2},
\end{equation}
where
\begin{equation}\label{mul res sig1}
	\gamma(p,\theta,\textbf{Q}):=\sup_{0\leq n' \leq n} \Big\{  \frac{2n'}{p} +\Big(\frac{2}{p}-(1-\theta)\Big)\mathfrak{d}_{d,n'}(\textbf{Q}) \Big\}.
\end{equation}
\end{lemma}

\begin{lemma}[\cite{GOZZK}, Lemma 5.5]\label{dec low var}
Let $2\leq q\leq p<\infty$. For every $d \geq 1$, $D>1$ and $\epsilon>0$, there exists $c=c(D,\epsilon)>0$ such that the following fact holds. For every sufficiently large $K$, there exist	
\begin{equation}\label{sec3 rel0}
	K^c \leq K_1 \leq K_2 \leq ...\leq K_D \leq K
\end{equation}
such that for every non-zero polynomial $P$ in $d$ variables of degree at most $D$, there exist collections of pairwise disjoint cubes $\mathcal{W}_j \subset \mathcal{P}(1/K_j)$, $j=1,2,...,D,$ such that 
\begin{equation}\label{sec3 rel1}
	N_{1/K}(Z_P) \cap [0,1]^d \subset  \bigcup_{j=1}^D \bigcup_{W \in \mathcal{W}_j} W
\end{equation}
and
\begin{equation}
	\Big\|\sum_{W \in \mathcal{W}_j} E^{\textbf{Q}}f_W\Big\|_{L^p(\mathbb{R}^{d+n})} \lesssim_{D,\textbf{Q},\epsilon,q,p} K_j^{\Gamma_{q,p}^{d-1}(\textbf{Q})+\epsilon} \Big(   \sum_{W \in \mathcal{W}_j}  \|E^{\textbf{Q}}f_W\|_{L^p(\mathbb{R}^{d+n})}^q\Big)^{\frac{1}{q}}.
\end{equation}
\end{lemma}

For $\textbf{Q}$ satisfying the strongly non-degenerate condition, we will first prove the following analogue of Proposition 3.1 in \cite{DZ}.

\begin{proposition}\label{sec2 le1}
Let $d,n \geq 1$, and $p=\frac{2(d-1+2n)}{d-1}$. Suppose that $\textbf{Q}$ is strongly non-degenerate. For any $0<\epsilon <1/100$, there are constant $C_\epsilon$ and $\delta=\epsilon^{100}$ such that the following fact holds for any $R \geq 1$ and every $f$ with ${\rm supp}f \subset B^d(0,1)$. Let $Y=\cup_{k=1}^M B_k$ be a union of lattice $K^2$-cubes in $B^{d+n}(0,R)$, where $K=R^\delta$. Suppose that 
	$$    \|E^{\textbf{Q}} f\|_{L^p(B_k)} \text{~is~essentially~a ~dyadic~constant~in~}   k=1,2,...,M.   $$
	Let $1\leq \alpha\leq d+n$ and
\begin{equation}\label{sign000}
\gamma:=\max_{\substack{  B^{d+n}(x',r)\subset B^{d+n}(0,R) \\ x'\in \mathbb{R}^{d+n},r\geq K^2 }} \frac{\#\{ B_k:B_k \subset B(x',r)\}}{r^\alpha}. 
\end{equation}
	Then 
	\begin{equation}\label{met1 eq1}
 \left\|  E^{\textbf{Q}} f \right\|_{L^p(Y)}   \leq C_\epsilon M^{-\frac{n}{d-1+2n}}\gamma^{\frac{n}{d-1+2n}}R^{\frac{(\alpha+n-1)n}{2(d-1+2n)}+\epsilon}\|f\|_{L^2}. 
	\end{equation}
\end{proposition}

\noindent \textit{Proof.} Let $K \in 2^\mathbb{N}$ be a dyadic integer. We decompose $B^{d}(0,1)$ into $K^{-1}$-cubes $\tau$. Then write $f=\sum_{\tau} f_\tau$, where $f_\tau=f \chi_\tau$. For a $K^2$-cube $B$ in $Y$, we define its significant set as
$$    \mathcal{S}(B) :=\Big\{ \tau: \| E^{\textbf{Q}} f _\tau\|_{L^p(B)} \geq \frac{1}{\# \tau} \max_{\tau'}\|  E^{\textbf{Q}} f_{\tau'}  \|_{L^p(B)} \Big\}.  $$
Running the algorithm from \cite[Propostion 5.6]{GOZZK}, together with Lemma \ref{dec low var}, we obtain: there exist $\mathcal{W}_{\iota,j}$ and $\mathcal{T}(B)$, such that each $\tau \in \mathcal{S}(B)$ can be covered by exactly one cube in
$$    \Big(  \bigcup_{\iota} \bigcup_{j=1}^d \mathcal{W}_{\iota,j}    \Big) \bigcup \mathcal{T}(B).  $$
Here each $\mathcal{W}_{\iota,j}$ denotes the collection of pairwise disjoint $1/K_j$-cubes near the zero set of one polynomial $P$ in $d$ variables of degree at most $d$, and $\mathcal{T}(B)$ is a $\theta$-uniform set of cubes. Besides this, $K_j$, $j=1,...,d$, satisfy (\ref{sec3 rel0}) and $ \iota \lesssim \log K$. Then by the  triangle inequality, one has
\begin{equation}\label{met1 eq2}
\|E^{\textbf{Q}}f\|_{L^p(B)} \leq \Big\|\sum_{\tau \in \mathcal{T}(B)} E^{\textbf{Q}}f_\tau\Big\|_{L^p(B)}+\Big\|  \sum_{\tau \in \cup_\iota \cup_j \mathcal{W}_{\iota,j}} E^{\textbf{Q}}f_\tau \Big\|_{L^p(B)} +\Big\|   \sum_{\tau \notin\mathcal{S}(B)  }  E^{\textbf{Q}}f_\tau\Big\|_{L^p(B)}.
\end{equation}
We say $B$ is broad if the first term dominates, otherwise we say $B$ is narrow. We denote the union of broad cubes $B $ in $ Y$ by $Y_{\text{broad}}$ and the union of narrow cubes $B $ in $ Y$ by $Y_{\text{narrow}}$. We call it the broad case if $Y_{\text{broad}}$ contains $\geq M/2$ many $K^2$-cubes, and the narrow case otherwise. 

\vskip0.5cm

\noindent \textbf{Broad case.} For each broad cube $B$, we have
\begin{equation*}
	\|E^{\textbf{Q}} f \|_{L^p(B)} \lesssim \Big\|\sum_{\tau \in \mathcal{T}(B)} E^{\textbf{Q}}f_\tau\Big\|_{L^p(B)} \leq K^d \max_{\tau \in \mathcal{T}(B)} \|E^{\textbf{Q}}f_\tau\|_{L^p(B)} \leq K^{2d} \prod_{j=1}^J \|E^{\textbf{Q}}f_{\tau_j}\|^{\frac{1}{J}}_{L^p(B)},
\end{equation*}
where $\{\tau_1,...,\tau_J\}$ is a $\theta$-uniform set. For a $K^2$-cube $B$ with center $x_B$, we decompose $B$ to balls of the form $B(x_B +v,2)$, where $v \in B(0,K^2)\cap \mathbb{Z}^{d+n}$. There is $v_j \in B(0,K^2) \cap \mathbb{Z}^{d+n}$ such that $\|E^{\textbf{Q}}f _{\tau_j}\|_{L^\infty(B)}$ is attained in $B(x_B+v_j,2)$. Define $v_j=(x_j,y_j)$ and 
$$  f_{\tau_j,v_j}(\xi):= f_{\tau_j}(\xi)  e^{2\pi i(x_j\cdot\xi+y_j\cdot \textbf{Q}(\xi))}. $$
By translation and locally constant property, one has 
\begin{align*}
	\|E^{\textbf{Q}} f \|_{L^p(B)} &  \lesssim K^{2d} \prod_{j=1}^J \|E^{\textbf{Q}}f_{\tau_j}\|^{\frac{1}{J}}_{L^p(B)} \leq K^{O(1)} \prod_{j=1}^J \|E^{\textbf{Q}}f_{\tau_j,v_j}\|^{\frac{1}{J}}_{L^p(B(x_B,2))} \\
	& \sim K^{O(1)} \left\|     \prod_{j=1}^{J} |E^{\textbf{Q}}f_{\tau_j,v_j}|^{\frac{1}{J}} \right\|_{L^p(B(x_B,2))} .
\end{align*}
Since there are only $K^{O(1)}$ choices for $\tau_j,v_j$ and $J$, by pigeonholing, there exist $\tilde{\tau_j},\tilde{v_j},\tilde{J}$ such that the above inequality holds for at least $\geq K^{-C}M$ many broad cubes. Now fix $\tilde{\tau_j},\tilde{v_j}$ and $\tilde{J}$.  We abbreviate $f_{\tilde{\tau_j},\tilde{v_j}}$ to $f_j$, and denote the collection of remaining broad cubes $B$ by $\mathcal{B}$. Next we sort $B \in \mathcal{B}$ by the value of $ \|\prod_{j=1}^{\tilde{J}} |E^{\textbf{Q}} f _j|^{\frac{1}{\tilde{J}}}\|_{L^\infty(B(x_B,2))}$: for dyadic number $A$, define
$$ \mathbb{Y}_A:=\left\{B \in \mathcal{B}: \left\|\prod_{j=1}^{\tilde{J}} |E^{\textbf{Q}} f _j|^{\frac{1}{\tilde{J}}}\right\|_{L^\infty(B(x_B,2))} \sim A\right\} .  $$
Let $Y_{A}$ be the union of the $K^2$-cubes $B$ in $\mathbb{Y}_{A}$. Without loss of generality, assume that $\|f\|_{L^2}=1$. We can further assume that $R^{-C} \leq A \leq 1$ for some constant $C$. So there are only $O(\log R) \leq O(K)$ choices on $A$. By dyadic pigeonholing, there exists a constant $\tilde{A}$ such that 
\begin{equation*}
\# \{  B: B\subset Y_{\tilde{A}}  \} \gtrsim (\log R)^{-1} \#  \mathcal{B}. 
\end{equation*}
Now we fix $\tilde{A}$, and denote $Y_{\tilde{A}}$ by $Y'$. We take $\theta \ll \epsilon$ and let 
\begin{equation}\label{sec2 sig1}
q=\max_{1\leq n' \leq n}
\Big\{2+\frac{2n'}{\mathfrak{d}_{d,n'}(\textbf{Q})}\Big\}
\end{equation}
denote the endpoint of (\ref{mul res sig1}) when $\theta=0$. Note that $2 \leq q \leq \frac{2(d+n)}{d} < p$ since $\textbf{Q}$ is strongly non-degenerate, it follows from Lemma \ref{mul th} that
\begin{align}
	\|E^{\textbf{Q}} f \|_{L^p(Y_{\text{broad}})} & \leq 	K^{O(1)}\|E^{\textbf{Q}} f \|_{L^p(Y')}    \nonumber \\ 
	& \leq K^{O(1)} \left\| \prod_{j=1}^{\tilde{J}} |E^{\textbf{Q}} f_j |^{\frac{1}{\tilde{J}}}  \right\|_{L^p(\cup_{B\subset Y'} B(x_B,2))}   \nonumber\\
	& \sim K^{O(1)} M^{\frac{1}{p}-\frac{1}{q}} \left\| \prod_{j=1}^{\tilde{J}} |E^{\textbf{Q}}f _j|^{\frac{1}{\tilde{J}}}  \right\|_{L^{q}(\cup_{B\subset Y'} B(x_B,2))}     \nonumber \\
	& \leq   K^{O(1)} M^{\frac{1}{p}-\frac{1}{q}}  \left\| \prod_{j=1}^{\tilde{J}} |E^{\textbf{Q}} f _j|^{\frac{1}{\tilde{J}}}  \right\|_{L^{q}(B_R)}   \nonumber\\
	& \lesssim    K^{O(1)} M^{\frac{1}{p}-\frac{1}{q}}  R^\epsilon\|f\|_{L^2}    \nonumber\\
	&  \leq K^{O(1)} M^{\frac{1}{p}-\frac{1}{2}}  \gamma^{\frac{1}{2}-\frac{1}{q}} R^{(\frac{1}{2}-\frac{1}{q})\alpha+\epsilon}\|f\|_{L^2} \label{s3equ2}\\
	&\lesssim K^{O(1)} M^{-\frac{n}{d-1+2n}}\gamma^{\frac{n}{d-1+2n}}R^{\frac{(\alpha+n-1)n}{2(d-1+2n)}+\epsilon}\|f\|_{L^2}.    \nonumber
\end{align}
Here we have used the facts $M\leq \gamma R^\alpha$ and $\gamma \geq K^{-2\alpha}$, which can be derived by the definition of $\gamma$.  

\vskip0.5cm

\noindent \textbf{Narrow case.}  For each narrow cube $B$, we have
\begin{equation*}
	\|E^{\textbf{Q}}f\|_{L^p(B)} \lesssim \Big\|  \sum_{\tau \in \cup_\iota \cup_j \mathcal{W}_{\iota,j}} E^{\textbf{Q}}f_\tau \Big\|_{L^p(B)} +\Big\|   \sum_{\tau \notin\mathcal{S}(B)  }  E^{\textbf{Q}}f_\tau\Big\|_{L^p(B)}.
\end{equation*}
If the second term dominates, by the definition of $\mathcal{S}(B)$, one gets
$$\Big\|   \sum_{\tau \notin\mathcal{S}(B)  }  E^{\textbf{Q}}f_\tau\Big\|_{L^p(B)}  \leq \sum_{\tau \notin\mathcal{S}(B)  }   \|   E^{\textbf{Q}}f_\tau\|_{L^p(B)} \leq \max_{\tau'}\|  E^{\textbf{Q}} f_{\tau'}  \|_{L^p(B)} .  $$
Therefore essentially only one $\tau$ plays a role in this case, and we can reduce it to the first term. From now on, we will only consider the first term. By the triangle inequality, one gets
\begin{equation}\label{add}
\|E^{\textbf{Q}}f\|_{L^p(B)} \lesssim \Big\|  \sum_{\tau \in \cup_\iota \cup_j W_{\iota,j}} E^{\textbf{Q}}f_\tau \Big\|_{L^p(B)} \leq  \sum_{\iota \lesssim \log K} \sum_{j=1}^{d}\Big\|  \sum_{\tau \in \mathcal{W}_{\iota,j}} E^{\textbf{Q}}f_\tau \Big\|_{L^p(B)}  .
\end{equation}

We plan to study (\ref{add}) through the same steps as in \cite[Section 3]{DZ}: do some dyadic pigeonholing, apply lower dimensional decoupling, and use rescaling and induction on scale. However, the lower dimensional decoupling Lemma \ref{dec low var} of the case can only reduce the support of $f$ to $K_j^{-1}$-scale, rather than $K^{-1}$-scale. In order to make the induction work, we need to consider the decomposition on $K_j^{-1}$-scale for the function in (\ref{add}), so that they satisfy the assumptions of Proposition \ref{sec2 le1} at $R/K_j^2$-scale after rescaling. Without loss of generality, we assume that the contribution from $\mathcal{W}_{1,1}$ dominates in (\ref{add}) for a fraction $\gtrsim (\log R)^{-1}$ of all $B $ in $Y_{\text{narrow}}$. We denote the union of such $B$ by $Y_1$. Recall that $\mathcal{W}_{1,1}$ is a collection of some pairwise disjoint $1/K_1$-cubes. On each $B $ in $ Y_1$, one has
\begin{equation}\label{s3equ1}
	\|E^{\textbf{Q}}f\|_{L^p(B)} \lesssim  \log K \Big\|  \sum_{\tau \in \mathcal{W}_{1,1}} E^{\textbf{Q}}f_\tau \Big\|_{L^p(B)}  \sim  \log K \Big\|  \sum_{W \in \mathcal{W}_{1,1}} E^{\textbf{Q}}f_W \Big\|_{L^p(B)}  ,
\end{equation}
where  each $W$ is a $1/K_1$-cube with center $c(W)$. We break $B^d(0,R)$ in the physical space into $R/K_1$-cubes $D$ with center $c(D)$. Write 
$$\sum_{W \in \mathcal{W}_{1,1}} f_W=\sum_{W,D} f_{W,D},$$ 
where each $f_{W,D}$ is of Fourier supported on $W \in \mathcal{W}_{1,1}$, and essentially supported on $D$. Then  
$$ \sum_{W \in \mathcal{W}_{1,1}} E^{\textbf{Q}}f_W=\sum_{W,D} E^{\textbf{Q}} f_{W,D}.$$ 
Through the stationary phase method, we can see that each $E^{\textbf{Q}} f_{W,D}$ restricted on $B_R$ is essentially supported on 
 $$\Box_{W,D}:=\Big\{    (x,y) \in B^{d+n}_R  : \Big|x+\sum_{j=1}^n y_j \nabla Q_j(c(W))-c(D)    \Big|\leq R/K_1 \Big\}.$$
Notice that, for each $W$, a given $B$ lies in exactly one box $\Box_{W,D}$, it follows from
 Lemma \ref{dec low var} and (\ref{sec3 rel0}) that
\begin{equation}\label{eeeeeeee}
\Big\|  \sum_{W \in \mathcal{W}_{1,1}} E^{\textbf{Q}}f_W \Big\|_{L^p(B)} \lesssim  K_1^{\epsilon^2} \Big(\sum_{ \Box}  \left\| E^{\textbf{Q}}f_\Box \right\|^2_{L^p(w_B)}\Big)^{\frac{1}{2}} , \quad\quad  \forall~ B \subset Y_1,
\end{equation}	
where $w_B$ is given by (\ref{sig1}), and $\Box:=\Box_{W,D}$.

Set $\tilde{R}=R/K_1^2$ and $\tilde{K}=\tilde{R}^\delta$. Tile $\Box$ by the
rectangle boxes $S$ with $d$ short sides of length $K_1\tilde{K}^2$ pointing in the same directions as the short sides of $\Box$, and $n$ long sides of length 
$K_1^2\tilde{K}^2$ pointing in complementary directions. We now perform dyadic pigeonholing argument to $S$ and $\Box$ which has similar spirit as in \cite[Section 3]{DZ}: 

\noindent (1) For each $\Box$, we sort $S \subset \Box$ that intersect $Y_{1}$ according to the value of $\|E^{\textbf{Q}} f_{\Box}\|_{L^p(S)}$ and the number of $K^2$-cubes in $Y_{1}$ contained in it: for dyadic numbers $\eta,\beta_1$, define
	$$    \mathbb{S}_{\Box,\eta,\beta_1}:=\Big\{ S\subset \Box:S~\text{contains}\sim\eta~\text{many}~K^2\text{-cubes~in~}Y_{1},~ \|E^{\textbf{Q}} f_{\Box}\|_{L^p(S)}\sim\beta_1  \Big\}.    $$
Let $Y_{\Box,\eta,\beta_1}$ be the union of the rectangular boxes $S$ in $\mathbb{S}_{\Box,\eta,\beta_1}$.  	
	
\noindent (2) For fixed $\eta$ and $\beta_1$, we sort $\Box$ according to the value of $\|f_\Box\|_{L^2}$ and the number $\#\mathbb{S}_{\Box,\eta,\beta_1}$: for dyadic numbers $\beta_2,M_1$, define
	$$   \mathbb{B}_{\eta,\beta_1,\beta_2,M_1}:=\Big\{ \Box: \|f_\Box\|_{L^2}\sim \beta_2,~ \# \mathbb{S}_{\Box,\eta,\beta_1}\sim M_1\Big\}.  $$

On each cube $B$ in $Y_{1}$, we have
$$\sum_{W \in \mathcal{W}_{1,1}} E^{\textbf{Q}}f_W =\sum_{\eta,\beta_1,\beta_2,M_1}\left( \sum_{\substack{\Box\in \mathbb{B}_{\eta,\beta_1,\beta_2,M_1}  \\ B\subset Y_{\Box,\eta,\beta_1}}} E^{\textbf{Q}} f_\Box  \right). $$
Under the hypothesis $\|f\|_{L^2}=1$, we can further assume that
\begin{align*}
	1\leq \eta \leq K^{O(1)}, \quad R^{-C} \leq \beta_1 \leq K^{O(1)}, \quad R^{-C} \leq \beta_2 \leq 1,    \quad 1 \leq M_1 \leq R^{O(1)}
\end{align*}
for some constant $C$. Thus there are only $O(\log R)$ choices for each dyadic number. By dyadic pigeonholing, these exist $\eta,\beta_1,\beta_2,M_1$ depending on $B$ such that (\ref{s3equ1}) becomes
\begin{equation}\label{absd}
	\|E^{\textbf{Q}}f\|_{L^p(B)} \lesssim  (\log R)^5     \left\| \sum_{\substack{\Box\in \mathbb{B}_{\eta,\beta_1,\beta_2,M_1}  \\ B\subset Y_{\Box,\eta,\beta_1}}} E^{\textbf{Q}}f_\Box  \cdot   \chi_{Y_{\Box,\eta,\beta_1}}  \right\|_{L^p(B)}  .
\end{equation}

Finally, we sort $B$ in $ Y_{1}$. Since there are only $O(\log R)$ choices on $\eta,\beta_1,\beta_2,M_1$, by pigeonholing, we can find $\tilde{\eta},\tilde{\beta}_1,\tilde{\beta}_2,\tilde{M}_1$ such that (\ref{absd}) holds for a fraction $\gtrsim (\log R)^{-4}$ of all $B$ in $Y_{1}$. We denote the union of such $B$ by $Y'$. From now on, we abbreviate $Y_{\Box,\tilde{\eta},\tilde{\beta_1}}$ and $\mathbb{B}_{\tilde{\eta},\tilde{\beta}_1,\tilde{\beta}_2,\tilde{M}_1}$ to $Y_\Box$ and $\mathbb{B}$, respectively. Next we further sort $B $ in $ Y'$ by the number $\# \{ \Box \in \mathbb{B}: B \subset Y_\Box  \}$: for dyadic number $\mu$, define
$$  \mathbb{Y}_\mu :=\Big\{  B \subset Y':  \#  \{ \Box \in \mathbb{B}: B\subset Y_\Box \}\sim\mu \Big\}.  $$
Let $Y_{\mu}$ be the union of $B$ in $\mathbb{Y}_\mu$. Using dyadic pigeonholing again, we can choose $\tilde{\mu}$ such that
\begin{equation}\label{flr eq7}
	\#\{ B : B \subset Y_{\tilde{\mu}} \} \gtrsim (\log R)^{-1} \#\{B : B \subset Y' \}. 
\end{equation}
From now on, we fix $\tilde{\mu}$, and denote $Y_{\tilde{\mu}}$ by $Y''$. 

For each cube $B$ in $Y''$, it follows from (\ref{eeeeeeee}) that
\begin{align*}
	\| E^{\textbf{Q}} f \|_{L^p(B)} \lesssim&  (\log R)^5 \Big\|  \sum_{\Box \in\mathbb{B} : B \subset Y_\Box}E^{\textbf{Q}} f_\Box \cdot \chi_{Y_{\Box}}  \Big\|_{L^p(B)}  \\
	\lesssim  & (\log R)^5  K_1^{\epsilon^2} \Big( \sum_{
		\Box \in \mathbb{B}: B \subset Y_\Box} \|E^{\textbf{Q}} f_\Box \|^2_{L^p(w_B)}  \Big)^{\frac{1}{2}}  \\
	\lesssim   &  (\log R)^5 K_1^{\epsilon^2}\tilde{\mu}^{\frac{1}{2}-\frac{1}{p}} \Big( \sum_{\Box \in \mathbb{B}: B \subset Y_\Box} \|E^{\textbf{Q}} f_\Box \|^p_{L^p(w_B)}  \Big)^{\frac{1}{p}}.    
\end{align*}
Summing $B$ in $ Y''$, we see that
\begin{align}
	\|E^{\textbf{Q}}f \|_{L^p(Y)} &\lesssim (\log R)^6 \|E^{\textbf{Q}} f \|_{L^p(Y'')}     \nonumber \\
	& \lesssim  (\log R)^{11} K_1^{\epsilon^2} \tilde{\mu}^{\frac{1}{2}-\frac{1}{p}} \Big(\sum_{\Box \in \mathbb{B}} \|E^{\textbf{Q}} f _\Box\|^p_{L^p(Y_\Box)}\Big)^{\frac{1}{p}} .   \label{11}
\end{align}
As for the term $ \|E^{\textbf{Q}}f _\Box\|_{L^p(Y_\Box)}$, using the change of of variables $$ \xi   \rightarrow c(\tau)+K_1^{-1}\xi,$$ 
we conclude
\begin{equation}\label{12}
 \|E^{\textbf{Q}} f _\Box\|_{L^p(Y_\Box)}=K_1^{-\frac{n}{d-1+2n}} \|E^{\textbf{Q}}g\|_{L^p(\tilde{Y})},
\end{equation}
for some $g$ with supp$g\subset B^d(0,1)$ and $\|g\|_{L^2}=\|f_\Box\|_{L^2}$, where $\tilde{Y}$ is the image of $Y_\Box$ under the new coordinates. We easily see that $\|E^{\textbf{Q}}g\|_{L^p(\tilde{Y})}$ just satisfies the conditions of Proposition \ref{sec2 le1} under the new scale $\tilde{R}$ and parameter $\tilde{M_1}$. Therefore we apply inductive hypothesis to obtain
\begin{equation}\label{13}
\|E^{\textbf{Q}}g\|_{L^p(\tilde{Y})} \lesssim  \tilde{M_1}^{-\frac{n}{d-1+2n}} \tilde{\gamma_1}^{\frac{n}{d-1+2n}} \tilde{R}^{\frac{(\alpha+n-1)n}{2(d-1+2n)}+\epsilon} \|g\|_{L^2},
\end{equation}
where $\tilde{\gamma_1}$ is given by
$$   \tilde{\gamma_1}=\max_{\substack{  B^{d+n}(x',r)\subset B^{d+n}(0,\tilde{R}) \\ x'\in \mathbb{R}^{d+n},r\geq \tilde{K}^2 }} \frac{\#\{ B_k:B_k \subset B(x',r)\}}{r^\alpha}.   $$
Combining (\ref{11}), (\ref{12}) with (\ref{13}), one has
\begin{align}
	\|E^{\textbf{Q}} f \|_{L^p(Y)} & \leq  K_1^{2\epsilon^2} K_1^{-\frac{n}{d-1+2n}} \tilde{\mu}^{\frac{n}{d-1+2n}}  \tilde{M_1}^{-\frac{n}{d-1+2n}} \tilde{\gamma_1}^{\frac{n}{d-1+2n}} \tilde{R}^{\frac{(\alpha+n-1)n}{2(d-1+2n)}+\epsilon} \Big(\sum_{\Box \in \mathbb{B}}  \|f_\Box\|_{L^2}^p \Big)^{\frac{1}{p}}\nonumber   \\
	& \lesssim   K_1^{2\epsilon^2}  K_1^{-\frac{n}{d-1+2n}} \Big(\frac{\tilde{\mu}}{\tilde{M_1} \#\mathbb{B}}\Big)^{\frac{n}{d-1+2n}}   \tilde{\gamma_1}^{\frac{n}{d-1+2n}} \tilde{R}^{\frac{(\alpha+n-1)n}{2(d-1+2n)}+\epsilon}  \|f\|_{L^2}. \label{14}
\end{align}
Here $\log R$ is absorbed in $K_1^{\epsilon^2}$ due to $ K^c \leq K_1 \leq K=R^\delta$.

We claim that the old parameters $M,\gamma$ and new parameters $\tilde{M_1},\tilde{\gamma_1}$ have the following relations:
\begin{equation}\label{15}
\frac{\tilde{\mu}}{\#\mathbb{B}}\lesssim \frac{(\log R)^6 \tilde{M_1}\tilde{\eta}}{M},\quad \quad\tilde{\gamma_1} \tilde{\eta}\lesssim  \gamma  K_1^{\alpha+n}.  
\end{equation}
The first estimate in (\ref{15}) can be derived through the same steps as in \cite[Section 3]{DZ}. Now we prove the second estimate:
\begin{align*}
\tilde{\gamma_1} \tilde{\eta}      &  \sim  \max_{  B_r\subset B_{\tilde{R}} : r\geq \tilde{K}^2 } \frac{\#\{ B_k:B_k \subset B(x',r)\}}{r^\alpha}  \cdot \#\{    B:B\subset S\cap Y_{1} \text{~for~any~fixed~}S\subset Y_\Box  \}   \\
 & \sim\max_{T_r\subset \Box:r\geq \tilde{K}^2} \frac{\#\{    
	S:S\subset Y_\Box \cap T_r\}}{r^\alpha}   \cdot \#\{    B:B\subset S\cap Y_{1} \text{~for~any~fixed~}S\subset Y_\Box  \}   \\
&\lesssim  \max_{T_r\subset \Box:r\geq \tilde{K}^2}  \frac{\# \{  B\subset Y : B \subset T_r\}}{r^\alpha} \leq \frac{K_1^n\gamma  (K_1r)^\alpha}{r^\alpha}=\gamma K_1^{\alpha+n}.
\end{align*}
Here in the second line, we change variables back, then each $B_r$ becomes $T_r$ which denotes the rectangular box in $\Box$ with $d$ short sides of length $K_1r$ pointing in the same directions as the short sides of $\Box$, and $n$ long sides of length $K_1^2 r$ pointing in complementary directions. In the last inequality, we cover $T_r$ by $\sim K_1^n$ finitely overlapping $K_1r$-balls. The proof of this claim is complete.

Finally, we put (\ref{15}) into (\ref{14}) to yield 
$$    \| E^{\textbf{Q}} f \|_{L^p(Y)}  \lesssim K_1^{2\epsilon^2-2\epsilon}M^{-\frac{n}{d-1+2n}}\gamma^{\frac{n}{d-1+2n}}R^{\frac{(\alpha+n-1)n}{2(d-1+2n)}+\epsilon}\|f\|_{L^2}.$$
Since $K^c \leq K_1 \leq K=R^\delta$, we can choose sufficiently large $R$, then the induction closes for the narrow case. 

\qed

Next we extend Proposition \ref{sec2 le1} to the case of general quadratic forms $\textbf{Q}$. Let $q$ be defined by (\ref{sec2 sig1}), and $p$ with $p \geq q$. We consider the following estimate:
$$   	\|E^{\textbf{Q}}f\|_{L^p(Y)}\lesssim M^{\frac{1}{p}-\frac{1}{2}} \gamma^{\frac{1}{2}-\frac{1}{p}}R^{w(p)+\epsilon}\|f\|_{L^2}.   $$
We repeat the proof of Proposition \ref{sec2 le1}: on the broad case, from (\ref{s3equ2}), we get 
$$   w(p) \geq \Big(\frac{1}{2}-\frac{1}{q}\Big)\alpha	;   $$
on the narrow case, to make the scale induction close, we need
 $$ w(p) \geq \frac{1}{2}\Big[\Gamma_p^{d-1}(\textbf{Q})+\frac{d+2n}{p}-\frac{d}{2}+(\alpha+n)\Big(\frac{1}{2}-\frac{1}{p}\Big)\Big].$$
Combining the constraints of both cases, we obtain the following results.
\begin{proposition}\label{sec2 pro2}
Let $d,n \geq 1$. Let $q$ be defined by (\ref{sec2 sig1}), and $p$ with $p \geq q$. Suppose $\textbf{Q}=(Q_1,...,Q_n)$ is a collection of quadratic forms in $d$ variables.  For any $0<\epsilon <1/100$, there are constant $C_\epsilon$ and $\delta=\epsilon^{100}$ such that the following fact holds for any $R \geq 1$ and every $f$ with ${\rm supp}f \subset B^d(0,1)$. Let  $Y=\cup_{k=1}^M B_k$ be a union of lattice $K^2$-cubes in $B^{d+n}(0,R)$, where $K=R^\delta$. Suppose that 
$$    \|E^{\textbf{Q}} f\|_{L^p(B_k)}  \text{~is~essentially~a ~dyadic~constant~in~}   k=1,2,...,M.   $$
Let $1\leq \alpha\leq d+n$ and 
$$  \gamma:=\max_{\substack{  B^{d+n}(x',r)\subset B^{d+n}(0,R) \\ x'\in \mathbb{R}^{d+n},r\geq K^2 }} \frac{\#\{ B_k:B_k \subset B(x',r)\}}{r^\alpha}.  $$
Then 
\begin{equation}
	\|E^{\textbf{Q}}f\|_{L^p(Y)}\leq C_\epsilon M^{\frac{1}{p}-\frac{1}{2}} \gamma^{\frac{1}{2}-\frac{1}{p}}R^{w(p)+\epsilon}\|f\|_{L^2},
\end{equation}
where $w(p)$ is given by 
\begin{equation}\label{12345}
   w(p)=\frac{1}{2}   \max   \Big\{   \max_{1\leq n' \leq n} \frac{\alpha n'}{\mathfrak{d}_{d,n'}(\textbf{Q})+n'},~     \Gamma_p^{d-1}(\textbf{Q})+\frac{d+n-\alpha}{p}+\frac{\alpha+n-d}{2}       \Big\}. 
\end{equation}   
\end{proposition}

Proposition \ref{sec2 pro2} implies the weighted restriction estimates associated with quadratic forms via the standard procedures  as in \cite[Section 2]{DZ}.

\begin{corollary}
	Let $d,n \geq 1$ and $p \geq 2$. Let  $\textbf{Q}=(Q_1,...,Q_n)$ be a sequence of quadratic forms defined on $\mathbb{R}^d$.  Let $\alpha \in (0,d+n]$, and $\mu$ be an $\alpha$-dimensional measure in $\mathbb{R}^{d+n}$. Suppose that $f$ satisfies ${\rm supp}f \subset B^d(0,1)$, then
	$$    \|E^{\textbf{Q}}f\|_{L^2(d\mu_R)} \leq C_\epsilon R^{w(p)+\epsilon}\|f\|_{L^2},   $$
	where $w(p)$ is given by 
$$   w(p)=\frac{1}{2}   \max   \Big\{   \max_{1\leq n' \leq n} \frac{\alpha n'}{\mathfrak{d}_{d,n'}(\textbf{Q})+n'},~     \Gamma_p^{d-1}(\textbf{Q})+\frac{d+n-\alpha}{p}+\frac{\alpha+n-d}{2}       \Big\}. $$
	In particular, if $\textbf{Q}$ is strongly non-degenerate, for $p=\frac{2(d-1+2n)}{d-1}$, we have 
	$$ 	\|E^{\textbf{Q}}f\|_{L^2(d\mu_R)} \leq C_\epsilon R^{\frac{(\alpha+n-1)n}{2(d-1+2n)}+\epsilon}\|f\|_{L^2}.   $$
\end{corollary}

In particular, when $d=3$ and $n=2$, using Corollary \ref{dqf3}, we obtain:

\begin{corollary}\label{cor1} 
Let $\alpha \in (0,5]$, and $\mu$ be an $\alpha$-dimensional measure in $\mathbb{R}^5$. Then
\begin{itemize}
\item[(1)] For $\textbf{Q}$ satisfying ${\rm (a)}$ and ${\rm (d)}$ in Theorem \ref{th1}, we have
\begin{equation}\label{c1ad}
	\|E^{\textbf{Q}}f\|_{L^2(d\mu_R)}\lesssim R^{\min\{   \frac{\alpha+1}{4},\frac{\alpha+2}{6} \}  +\epsilon}\|f\|_{L^2}.   
\end{equation}
\item[(2)] For $\textbf{Q}$ satisfying ${\rm (b)}$ in Theorem \ref{th1}, we have
\begin{equation}\label{c1b}
	\|E^{\textbf{Q}}f\|_{L^2(d\mu_R)}\lesssim R^{\max\left\{ \frac{\alpha}{4},  \min\{    \frac{\alpha+1}{4}, \frac{\alpha+3}{8} \}    \right\}+\epsilon}\|f\|_{L^2}.  
\end{equation}
\item[(3)] For $\textbf{Q}$ satisfying ${\rm (c)}$ in Theorem \ref{th1}, we have
\begin{equation}\label{c1c}
	\|E^{\textbf{Q}}f\|_{L^2(d\mu_R)}\lesssim R^{ \frac{\alpha+1}{4} +\epsilon}\|f\|_{L^2}.  
\end{equation}
\item[(4)] Suppose that $\textbf{Q}$ satisfies ${\rm (e)}$ in Theorem \ref{th1}. If $\textbf{Q}$ satisfies  $\mathfrak{d}_{2,1}(\textbf{Q})=\mathfrak{d}_{1,2}(\textbf{Q})=1$ which is strongly non-degenerate, then
\begin{equation}\label{c1e}
	\|E^{\textbf{Q}}f\|_{L^2(d\mu_R)}\lesssim R^{\frac{\alpha+1}{6}+\epsilon}\|f\|_{L^2}.
\end{equation}
If $\textbf{Q}$ satisfies  $\mathfrak{d}_{2,1}(\textbf{Q})=1$ and  $\mathfrak{d}_{1,2}(\textbf{Q})=0$, then
\begin{equation}
	\|E^{\textbf{Q}}f\|_{L^2(d\mu_R)}\lesssim R^{\frac{\alpha+1}{5}  +\epsilon}\|f\|_{L^2}.     
\end{equation}
If $\textbf{Q}$ satisfies  $\mathfrak{d}_{2,1}(\textbf{Q})=0$ and  $\mathfrak{d}_{1,2}(\textbf{Q})=1$,  then
\begin{equation}
	\|E^{\textbf{Q}}f\|_{L^2(d\mu_R)}\lesssim R^{ \min\{ \frac{\alpha+1}{4},   \frac{\alpha+3}{8} \}   +\epsilon}\|f\|_{L^2}.
\end{equation}
If $\textbf{Q}$ satisfies  $\mathfrak{d}_{2,1}(\textbf{Q})=\mathfrak{d}_{1,2}(\textbf{Q})=0$, then
\begin{equation}
	\|E^{\textbf{Q}}f\|_{L^2(d\mu_R)}\lesssim R^{\min\{ \frac{\alpha+1}{4}, \frac{\alpha+2}{6} \} +\epsilon}\|f\|_{L^2}.       
\end{equation}
\end{itemize}
\end{corollary}

\section{A direct application of the $\ell^2$ decoupling}

In the previous section, we have used the Du-Zhang method to obtain a few weighted restriction estimates associated with quadratic forms. In this section, instead of using the broad-narrow analysis, we will use the $\ell^2$ decoupling for quadratic forms directly. For the standard parabolic case as in \cite{DZ}, such treatment can not lead to better results than that from the previous section. The main reasons come from its nice multilinear restriction estimate and lower dimensional decoupling. However,  for some degenerate quadratic forms, the situations are quite different. For example, from Corollary \ref{dqf3}, we see that $\textbf{Q}$ satisfying (a) and (d) in Theorem \ref{th1} have the same decoupling endpoint exponents both in original and lower dimensions.

\begin{proposition}\label{sec3 pro1}
	Let $d,n \geq 1$ and $p \geq 2$. Suppose that $\textbf{Q}=(Q_1,...,Q_n)$ is a collection of quadratic forms in $d$ variables. For any $0<\epsilon <1/100$, there are constant $C_\epsilon$ and $\delta=\epsilon^{100}$ such that the following fact holds for any $R \geq 1$ and every $f$ with ${\rm supp}f \subset B^d(0,1)$. Let $Y=\cup_{k=1}^M B_k$ be a union of lattice $K^2$-cubes in $B^{d+n}(0,R)$, where $K=R^\delta$. Suppose that 
	$$    \|E^{\textbf{Q}} f\|_{L^p(B_k)}  \text{~is~essentially~a ~dyadic~constant~in~}   k=1,2,...,M.   $$
	Let $1\leq \alpha\leq d+n$ and 
	$$  \gamma:=\max_{\substack{  B^{d+n}(x',r)\subset B^{d+n}(0,R) \\ x'\in \mathbb{R}^{d+n},r\geq K^2 }} \frac{\#\{ B_k:B_k \subset B(x',r)\}}{r^\alpha}.  $$
	Then 
	\begin{equation}
		\|E^{\textbf{Q}}f\|_{L^p(Y)}\leq C_\epsilon M^{\frac{1}{p}-\frac{1}{2}}\gamma^{\frac{1}{2}-\frac{1}{p}}R^{w(p)+\epsilon}\|f\|_{L^2},
	\end{equation}
	where $w(p)$ is given by 
\begin{equation}\label{123456}
  w(p)=\frac{1}{2}   \Big[\Gamma_p^{d}(\textbf{Q})+\frac{d+n-\alpha}{p}+\frac{\alpha+n-d}{2}    \Big]       . 
\end{equation}	
In particular, if $\textbf{Q}$ is strongly non-degenerate, for $p=\frac{2(d+2n)}{d}$, we have 
$$ 		\|E^{\textbf{Q}}f\|_{L^p(Y)}\leq C_\epsilon M^{-\frac{n}{d+2n}} \gamma^{\frac{n}{d+2n}}R^{\frac{(\alpha+n)n}{2(d+2n)}+\epsilon}\|f\|_{L^2}.   $$	
\end{proposition}

We can prove this proposition by repeating the argument of  Proposition \ref{sec2 le1} in the narrow case, together with the $\ell^2$ decoupling for quadratic forms. Then similar to the argument as in the previous section, we can use Proposition \ref{sec3 pro1} to get a few weighted restriction estimates associated with quadratic forms.

\begin{corollary}\label{buchong2}
	Let $d,n \geq 1$ and $p \geq 2$. Let  $\textbf{Q}=(Q_1,...,Q_n)$ be a sequence of quadratic forms defined on $\mathbb{R}^d$.  Let $\alpha \in (0,d+n]$, and $\mu$ be an $\alpha$-dimensional measure in $\mathbb{R}^{d+n}$. Suppose that $f$ satisfies ${\rm supp}f \subset B^d(0,1)$, then
	$$    \|E^{\textbf{Q}}f\|_{L^2(d\mu_R)} \leq C_\epsilon R^{w(p)+\epsilon}\|f\|_{L^2},   $$
	where $w(p)$ is given by 
   $$w(p)=\frac{1}{2}   \Big[\Gamma_p^{d}(\textbf{Q})+\frac{d+n-\alpha}{p}+\frac{\alpha+n-d}{2}    \Big]       . $$
	In particular, if $\textbf{Q}$ is strongly non-degenerate, for $p=\frac{2(d+2n)}{d}$, we have 
	$$ 	\|E^{\textbf{Q}}f\|_{L^2(d\mu_R)} \leq C_\epsilon R^{\frac{(\alpha+n)n}{2(d+2n)}+\epsilon}\|f\|_{L^2}.   $$
\end{corollary}

In particular, when $d=3$ and $n=2$, using Corollary \ref{dqf3}, we obtain:
\begin{corollary}\label{cor2} 
Let $\alpha \in (0,5]$, and $\mu$ be an $\alpha$-dimensional measure in $\mathbb{R}^5$.	Then
	\begin{itemize}
		\item[(1)] For $\textbf{Q}$ satisfying ${\rm (a)}$ and ${\rm (d)}$ in Theorem \ref{th1}, we have
		\begin{equation}\label{c2ad}
			\|E^{\textbf{Q}}f\|_{L^2(d\mu_R)}\lesssim R^{\frac{\alpha+2}{6}  +\epsilon}\|f\|_{L^2}.   
		\end{equation}
		\item[(2)] For $\textbf{Q}$ satisfying ${\rm (b)}$ in Theorem \ref{th1}, we have
		\begin{equation}\label{c2b}
			\|E^{\textbf{Q}}f\|_{L^2(d\mu_R)}\lesssim R^{\min\left\{ \frac{\alpha+2}{6},\frac{\alpha+3}{8}   \right\}+\epsilon}\|f\|_{L^2}.  
		\end{equation}
			\item[(3)] For $\textbf{Q}$ satisfying ${\rm (c)}$ in Theorem \ref{th1}, we have
	\begin{equation}\label{c2c}
		\|E^{\textbf{Q}}f\|_{L^2(d\mu_R)}\lesssim R^{\frac{\alpha+2}{5}+\epsilon}\|f\|_{L^2}.  
	\end{equation}
		\item[(4)] Suppose that $\textbf{Q}$ satisfies ${\rm (e)}$ in Theorem \ref{th1}. If $\textbf{Q}$ satisfies  $\mathfrak{d}_{2,1}(\textbf{Q})=\mathfrak{d}_{1,2}(\textbf{Q})=1$ which is strongly non-degenerate, then
		\begin{equation}\label{c2e}
			\|E^{\textbf{Q}}f\|_{L^2(d\mu_R)}\lesssim R^{\frac{\alpha+2}{7}+\epsilon}\|f\|_{L^2}.
		\end{equation}
		If $\textbf{Q}$ satisfies  $\mathfrak{d}_{2,1}(\textbf{Q})=1$ and  $\mathfrak{d}_{1,2}(\textbf{Q})=0$, then
		\begin{equation}
	\|E^{\textbf{Q}}f\|_{L^2(d\mu_R)}\lesssim R^{\frac{\alpha+2}{6}  +\epsilon}\|f\|_{L^2}.   
\end{equation}
		If $\textbf{Q}$ satisfies  $\mathfrak{d}_{2,1}(\textbf{Q})=0$ and  $\mathfrak{d}_{1,2}(\textbf{Q})=1$, then
		\begin{equation}
	\|E^{\textbf{Q}}f\|_{L^2(d\mu_R)}\lesssim R^{\min\left\{ \frac{\alpha+2}{6},\frac{\alpha+3}{8}   \right\}+\epsilon}\|f\|_{L^2}.  
\end{equation}
		If $\textbf{Q}$ satisfies  $\mathfrak{d}_{2,1}(\textbf{Q})=\mathfrak{d}_{1,2}(\textbf{Q})=0$, then
		\begin{equation}
			\|E^{\textbf{Q}}f\|_{L^2(d\mu_R)}\lesssim R^{\frac{\alpha+2}{6}  +\epsilon}\|f\|_{L^2}.      
		\end{equation}
	\end{itemize}
\end{corollary}

Now we compare Corollary \ref{cor1} and Corollary \ref{cor2}. For the degenerate case, the results from Corollary \ref{cor2} may be better than that from Corollary \ref{cor1} though the argument of this section seems rough. For example, for $\textbf{Q}$ satisfying (b) in Theorem \ref{th1}, the power $\alpha/4$ in $(\ref{c1b})$ comes from the corresponding multilinear restriction estimates, which reveals that the broad-narrow argument can not help us to gain better results in this case.
Nonetheless, for the non-degenerate case, the estimates in Corollary \ref{cor2} are weaker indeed than that in Corollary \ref{cor1} such as $\textbf{Q}$ satisfying the strongly non-degenerate condition.

Finally, we point out that Corollary \ref{buchong2} can be also derived by the refined Strichartz estimate. We leave the details in Appendix 1.

\section{$k$-linear restriction estimates for quadratic forms}

Our goal in this section is to build one class of $k$-linear restriction estimates for quadratic forms. A central tool is the Brascamp-Lieb inequality. Such inequality gives a general form for a host of important inequalities such as H\"older, Young, and Loomis-Whitney. There are several versions on such inequality, including the scale-invariant, discrete, and local Brascamp-Lieb inequalities, which have been characterized in \cite{BCCT08,BCCT10}. Recently, the scale-dependent version of the Brascamp-Lieb inequality, which unifies the scale-invariant, discrete, and local Brascamp-Lieb inequalities, was obtained in  \cite{Mal22}. It was used to obtain the sharp decoupling inequalities for quadratic forms in \cite{GOZZK} for the first time.

\begin{definition}\label{bl def}
Let $m,m'\in \mathbb{N}$. Let $(V_j)^M_{j=1}$ be a tuple of linear subspaces $V_j \subset \mathbb{R}^m$ of dimension $m'$. For a linear subspace $V \subset \mathbb{R}^m$, let $\pi_V: \mathbb{R}^m \rightarrow V$ denote the orthogonal projection onto $V$. For $0\leq \alpha \leq M$ and $R \geq 1$, we denote by $\rm{BL}$$((V_j)^M_{j=1}, \alpha, R, \mathbb{R}^m)$ the smallest constant such that the inequality
\begin{equation}
\int_{[-R,R]^m} \prod_{j=1}^{M} f_j(\pi_{V_j}(x))^{\frac{\alpha}{M}} dx \leq {\rm BL}((V_j)^M_{j=1}, \alpha, R, \mathbb{R}^m) \prod_{j=1}^{M} \Big(   \int_{V_j} f_j(x_j) dx_j \Big)^{\frac{\alpha}{M}}
\end{equation} 
holds for any functions $f_j: V_j \rightarrow [0,\infty)$ that are constant at scale 1, in the sense that $V_j$ can be partitioned into cubes with unit side length on each of which $f_j$ is constant. We abbreviate $\rm{BL}$$((V_j)^M_{j=1}, \alpha, R):=$$\rm{BL}$$((V_j)^M_{j=1}, \alpha, R,\mathbb{R}^m)$ when $\mathbb{R}^m$ is clear in the context.
\end{definition}

\begin{theorem}[\cite{Mal22}, Theorem 3]\label{bl th} 
In the notation of Definition \ref{bl def}, fix a tuple $(V_j)_{j=1}^M$ and an exponent $1 \leq \alpha \leq M$. Let 
$$ \kappa:=\sup_{V \subset \mathbb{R}^m} \Big( \dim V -\frac{\alpha}{M} \sum_{j=1}^{M} \dim \pi_{V_j} V  \Big),  $$	
where the supremum is taken over all linear subspaces of $\mathbb{R}^m$.
Then there exists a constant $C_0 <\infty$ and a neighborhood of the tuple $(V_j)_{j=1}^M$ in the $M$-th power of the Grassmanian manifold of all linear subspaces of dimension $m'$ of $\mathbb{R}^m$ such that, for any tuple $(\tilde{V}_j)_{j=1}^M$ in this neighborhood and any $R \geq 1$, we have
\begin{equation}
{\rm BL}((\tilde{V}_j)^M_{j=1}, \alpha, R) \leq C_0 R^\kappa.
\end{equation}
\end{theorem}

Now we introduce a $k$-transverse condition on quadratic forms. Let  $\textbf{Q}=(Q_1,...,Q_n)$ be a sequence of quadratic forms defined on $\mathbb{R}^d$. Recall that $S_{\textbf{Q}}=\{   (\xi,\textbf{Q}(\xi) ):\xi \in [0,1]^d \}$ denotes the associated graph of $\textbf{Q}$. The normal space at the point $\xi$ of $S_{\textbf{Q}}$ is given by
\begin{equation}\label{bcbc1}
N_\xi(\textbf{Q}):= N_\xi:=\text{span} \{  (\nabla Q_j(\xi),-e_j) :j=1,...,n \}:=\text{span} \{ \textbf{n}_\xi^j :j=1,...,n \}, 
\end{equation} 
and the tangent space at the point $\xi$ of $S_{\textbf{Q}}$ is given by
\begin{equation}\label{bcbc2}
V_\xi (\textbf{Q}):=V_\xi:=\text{span} \{  (e_i, \partial_i \textbf{Q}_i(\xi)) :i=1,...,d \}. 
\end{equation}
\begin{definition}\label{k trans def}
Let $k \geq 2$. Let  $\textbf{Q}=(Q_1,...,Q_n)$ be a sequence of quadratic forms defined on $\mathbb{R}^d$. Given subsets $U_1, ..., U_k$ in $[0,1]^d$, they are called $k$-transverse if for any $\gamma_j \in U_j$, $j=1, ..., k,$ we have
\begin{equation}\label{k trans con}
		\Big| \mathbf{n}^1_{\gamma_1} \wedge...\wedge \mathbf{n}^n_{\gamma_1}\wedge...\wedge \mathbf{n}^1_{\gamma_k} \wedge...\wedge \mathbf{n}^n_{\gamma_k} 	\Big| \gtrsim 1.
\end{equation}
\end{definition}

\begin{remark}
The possible values of $k$ depend on the values of $d$,  $n$, and the characteristics of $\textbf{Q}$. In general, if $\textbf{Q}$ satisfies $\mathfrak{d}_{d,1}(\textbf{Q})>0$, then $k$ can be taken in the range $2\leq k \leq \frac{\mathfrak{d}_{d,n}(\textbf{Q})+n}{n} $. In particular, if $\mathfrak{d}_{d,n}(\textbf{Q})<n$, then	(\ref{k trans con}) will not happen for any $k$.
\end{remark}

Under this transverse condition, we can use Theorem \ref{bl th} to gain the following $k$-linear Kakeya estimates for quadratic forms.

\begin{theorem}[$k$-linear Kakeya estimates for quadratic forms]\label{k linear kakeya}
	Let $d,n \geq 1$, $k\geq 2$ and $R \geq 1$. Let $\mathbb{T}_1,...,\mathbb{T}_k$ be the collections of rectangular boxes with $d$ short sides of length $1$ and $n$ long sides of length $R$. If for any $T_j \in \mathbb{T}_j$, $j=1, ..., k$, the directions of $n$ long sides of $T_j$ satisfy (\ref{k trans con}), then
	\begin{equation}\label{s5kakeya}
		\left\|  \prod_{j=1}^{k} \Big(  \sum_{T_j \in \mathbb{T}_j} \chi_{T_j} \Big)^{\frac{1}{k}} \right\|_{L^{\frac{k}{k-1}}(\mathbb{R}^{d+n})} \lesssim R^\epsilon \prod_{j=1}^{k} (\# \mathbb{T}_j)^{\frac{1}{k}}.	
	\end{equation} 
\end{theorem}

\noindent \textit{Proof.}  We apply Theorem \ref{bl th} to the case: $m=d+n$, $M=k$, and  each $V_j$, $j=1,...,k$, is an $n$-dimensional subspace which is consisted of the directions of $n$ long sides of a fixed  $T_j\in \mathbb{T}_j$. Then
\begin{equation}\label{k rans eq0}
	\kappa =\sup_{V\subset \mathbb{R}^{d+n}} \Big[
	\dim V -\frac{\alpha}{k}\Big(\sum_{j=1}^k  \dim\pi_{V_j}V\Big)\Big].
\end{equation} 
Now we consider the linear map $\pi_{V_j}: V \rightarrow V_j $, which is the
orthogonal projection onto $W$ restricted to $V$. Using the rank-nullity theorem, we get 
$$   \dim V = \dim\pi_{V_j}V + \dim (V \cap {V_j}^\perp) .$$
Thus
\begin{align*}
\dim V -\frac{\alpha}{k}\Big(\sum_{j=1}^k  \dim\pi_{V_j}V\Big) &= \dim V - \frac{\alpha}{k}\Big[   k \dim V -\sum_{j=1}^{k}\dim (V \cap {V_j}^\perp)   \Big].
\end{align*}
On the other hand, it follows from the dimensional formula and transverse condition (\ref{k trans con}) that
\begin{align*}
	\dim (V \cap V_1^\perp) +\dim (V \cap V_2^\perp)&= \dim (V \cap (V_1^\perp+V_2^\perp)) +\dim(V \cap V_1^\perp \cap V_2^\perp)   \\
&	= \dim (V \cap (V_1^\perp+V_2^\perp)).
\end{align*}
Using it repetitively, we obtain
$$ \sum_{j=1}^{k}\dim (V \cap V_j^\perp)  = \dim (V \cap (V_1^\perp+...+V_k^\perp))\leq \dim V.   $$
It follows that
$$  \kappa \leq \sup_{V\subset \mathbb{R}^{d+n}} \Big(1-\frac{k-1}{k} \cdot \alpha\Big)\dim V.   $$
Choosing $\alpha=\frac{k}{k-1}$, we establish the Brascamp-Lieb inequality under the transverse condition (\ref{k trans con}):
 \begin{equation}\label{bl eq3}
\int_{\mathbb{R}^{d+n}}   \prod_{j=1}^{k} f_j(\pi_{V_j}(x))^{\frac{1}{k-1}}dx \lesssim \prod_{j=1}^{k} \|f_j\|_{L^1}^{\frac{k}{k-1}}.  
 \end{equation}
Then (\ref{s5kakeya}) can be derived immediately by following the standard procedures of the Loomis-Whitney inequality implying the  multilinear Kakeya estimate as in \cite{Gut15}. 

\qed

Finally, via the argument of passing from the multilinear Kakeya estimate to the multilinear restriction estimate as in \cite{BCT06}, we obtain the following $k$-linear restriction estimates for quadratic forms, which can be seen as a generalization of the case of co-dimension 1.

\begin{theorem}[$k$-linear restriction estimates for quadratic forms]\label{klin res}
	Let $k \geq 2$. Let $\textbf{Q}=(Q_1,...,Q_n)$ be a sequence of quadratic forms defined on $\mathbb{R}^d$ with $d \geq 1$ and $n \geq 1$. Suppose that $\tau_1, ..., \tau_k$ are $k$-transverse in $[0,1]^d$. Then for each $f_j:\tau_j \rightarrow \mathbb{C}$, and any $R \geq 1$, we have
	\begin{equation}\label{klin res eq}
		\left\|   \prod_{j=1}^{k}  |E^{\textbf{Q}}f_j|^{\frac{1}{k}}   \right\|_{L^{\frac{2k}{k-1}}(B_R)} \lesssim R^\epsilon\prod_{j=1}^{k} \|f_j\|_{L^2}^{\frac{1}{k}}. 
	\end{equation}
\end{theorem}

\begin{remark}
When $k=2$, several special cases of Theorem \ref{klin res}, such as $d=3,~n=2,$ and $d=n=3$, have been derived by Guo-Oh \cite{GO} and Oh \cite{O18}. They used $L^4$ bi-orthogonal method to build the bilinear restriction estimates without $R^\epsilon$ loss.
\end{remark}

\section{The Du-Zhang method with the dimension 2 as a division}

Under the mechanism of the Du-Zhang method, we used the multilinear restriction estimates and lower dimensional $\ell^2$ decoupling to attain Proposition \ref{sec2 pro2} in Section 3, and the original dimensional  $\ell^2$ decoupling to attain Proposition \ref{sec3 pro1} in Section 4, respectively. The advantage of both methods is that they can be used widely for quadratic forms with any $d\geq 1$ and $n\geq 1$. However, note that (\ref{12345}) and (\ref{123456}) only depend on $\mathfrak{d}_{d',n'}(\textbf{Q})$, the two methods can't make a distinction for two quadratic forms with the same values of $\mathfrak{d}_{d',n'}(\textbf{Q})$ for all $0 \leq d' \leq d$ and all $0\leq n'\leq n$, such as $(\xi_1\xi_2,\xi_2^2+\xi_1\xi_3)$ and $ (\xi_1\xi_2,\xi_1^2\pm\xi^2_3)$ in Theorem \ref{th1}. In this section, under the framework of
 the broad-narrow analysis with the dimension 2 as a division, we use the bilinear restriction estimates built in the previous section in the broad case, and the adaptive $\ell^2$ decoupling in the narrow case. This argument can help us to make full use of the special degenerate characteristic of every quadratic form.

In the remainder of this section, we consider $\textbf{Q}$ in Theorem \ref{th1}.
We will investigate the following estimate:
\begin{equation}\label{aim}
	\|E^{\textbf{Q}}f\|_{L^p(Y)}\lesssim M^{\frac{1}{p}-\frac{1}{2}} \gamma^{\frac{1}{2}-\frac{1}{p}}R^{w(p)+\epsilon}\|f\|_{L^2},
\end{equation}
where $Y$, $M$ and $\gamma$ are defined as in Proposition \ref{sec2 pro2}.

We run the broad-narrow analysis. Divide $B^3(0,1)$ into $K^{-1}$-cubes $\tau$, and  write $f=\sum_\tau f_\tau$, where $f_\tau=f\chi_\tau$. For each $K^2$-cube $B$ in $Y$, define its  significant set as
$$    \mathcal{S}(B) :=\Big\{ \tau: \| E^{\textbf{Q}} f _\tau\|_{L^p(B)} \geq \frac{1}{100\# \tau} \|  E^{\textbf{Q}} f  \|_{L^p(B)} \Big\}.  $$
Then
$$    \Big\|  \sum_{\tau \in \mathcal{S}(B)} E^{\textbf{Q}}f_\tau \Big\|_{L^p(B)} \sim \|E^{\textbf{Q}}f\|_{L^p(B)}.  $$
We say $B$ is broad if there exist $\tau_1,\tau_2 \in \mathcal{S}(B)$ such that for any $\xi \in \tau_1$, $\eta \in \tau_2$
\begin{equation}\label{s6e1}
		\Big| \textbf{n}^1_{\xi} \wedge \textbf{n}^2_{\xi}\wedge \textbf{n}^1_{\eta} \wedge \textbf{n}^2_{\eta}  	\Big| \gtrsim K^{-O(1)}.
\end{equation}
Otherwise, we say $B$ is narrow. We denote the union of broad cubes $B $ in $Y$ by $Y_{\text{broad}}$ and the union of narrow cubes $B $ in $Y$ by $Y_{\text{narrow}}$. We call it the broad case if $Y_{\text{broad}}$ contains $\geq M/2$ many $K^2$-cubes, and the narrow case otherwise. 

\vskip0.5cm

\noindent \textbf{Broad case.} For each broad cube $B$, there exist $\tau_1,\tau_2 \in \mathcal{S}(B) $ satisfying the transverse condition (\ref{s6e1}) such that
\begin{equation*}
\int_B |E^{\textbf{Q}}f|^p \lesssim K^{3p}  \prod_{j=1}^{2} \Big( \int_B |E^{\textbf{Q}}f_{\tau_j}|^p  \Big)^{\frac{1}{2}}.
\end{equation*}
Repeating the argument as in Section 3, together with locally constant property, one has
\begin{equation*}
\int_B |E^{\textbf{Q}}f|^p \lesssim K^{O(1)}   \int_{B(x_B,2)} 	 \prod_{j=1}^{2} |E^{\textbf{Q}}f_{\tau_j,v_j}|^{\frac{p}{2}} ,
\end{equation*}
for some $v_j \in B(0,K^2) \cap \mathbb{Z}^{5}$.
Since there are only $K^{O(1)}$ choices for $\tau_j$ and $v_j$, by pigeonholing, there exist $\tilde{\tau_j}$ and $\tilde{v_j}$ such that the above inequality holds for at least $\geq K^{-C}M$ many broad balls. Now fix $\tilde{\tau_j}$ and $\tilde{v_j}$.  We abbreviate $f_{\tilde{\tau_j},\tilde{v_j}}$ to $f_j$, and denote the collection of remaining broad balls $B$ by $\mathcal{B}$. Next we sort $B \in \mathcal{B}$ by the value of $ \|\prod_{j=1}^{2} |E^{\textbf{Q}} f _j|^{\frac{1}{2}}\|_{L^\infty(B(x_B,2))}$: for dyadic number $A$, define
$$ \mathbb{Y}_A:=\left\{B \in \mathcal{B}: \left\|\prod_{j=1}^{2} |E^{\textbf{Q}} f _j|^{\frac{1}{2}}\right\|_{L^\infty(B(x_B,2))} \sim A\right\} .  $$
Let $Y_{A}$ be the union of the $K^2$-cubes $B$ in $\mathbb{Y}_{A}$. Without loss of generality, assume that $\|f\|_{L^2}=1$. We can further assume that $R^{-C} \leq A \leq 1$ for some constant $C$. So there are only $O(\log R) \leq O(K)$ choices on $A$. By dyadic pigeonholing, there exists a constant $\tilde{A}$ such that 
\begin{equation*}
	\# \{  B: B\subset Y_{\tilde{A}}  \} \gtrsim (\log R)^{-1}\#  \mathcal{B}. 
\end{equation*}
Now we fix $\tilde{A}$, and denote $Y_{\tilde{A}}$ by $Y'$. When $p \geq 4$, it follows from Theorem \ref{klin res} that
\begin{align*}
	\|E^{\textbf{Q}} f \|_{L^p(Y_{\text{broad}})} & \leq 	K^{O(1)}\|E^{\textbf{Q}} f \|_{L^p(Y')}   \\ 
	& \leq K^{O(1)} \left\| \prod_{j=1}^{2} |E^{\textbf{Q}} f_j |^{\frac{1}{2}}  \right\|_{L^p(\cup_{B\subset Y'} B(x_B,2))} \\
	& \sim K^{O(1)} M^{\frac{1}{p}-\frac{1}{4}} \left\| \prod_{j=1}^{2} |E^{\textbf{Q}}f _j|^{\frac{1}{2}}  \right\|_{L^{4}(\cup_{B\subset Y'} B(x_B,2))}    \\
	& \leq   K^{O(1)} M^{\frac{1}{p}-\frac{1}{4}}  \left\| \prod_{j=1}^{2} |E^{\textbf{Q}} f _j|^{\frac{1}{2}}  \right\|_{L^{4}(B_R)} \\
	& \lesssim    K^{O(1)} M^{\frac{1}{p}-\frac{1}{4}} R^\epsilon \|f\|_{L^2}  \\
	& =K^{O(1)}  M^{\frac{1}{p}-\frac{1}{2}}  M^{\frac{1}{4}}  R^\epsilon   \|f\|_{L^2}   \\
	&  \leq  K^{O(1)} M^{\frac{1}{p}-\frac{1}{2}}  \gamma^{\frac{1}{4}} R^{\frac{\alpha}{4}+\epsilon}\|f\|_{L^2} \\
	&\lesssim  K^{O(1)}  M^{\frac{1}{p}-\frac{1}{2}}  \gamma^{\frac{1}{2}-\frac{1}{p}} R^{\frac{\alpha}{4}+\epsilon}\|f\|_{L^2}.  
\end{align*}
The argument of the broad case is fitted for all $\textbf{Q}$ in Theorem \ref{th1}.

\vskip0.5cm

\noindent \textbf{Narrow case.} For a fixed narrow cube $B$, it is worth noting that different $\textbf{Q}$ possess different characteristics on $\mathcal{S}(B)$. We are going to study $\textbf{Q}$ in Theorem \ref{th1} in turn.

We start with $\textbf{Q}$ satisfying (a) in Theorem \ref{th1}. It suffices to consider $\textbf{Q}= (\xi_1^2,\xi_2^2+\xi_1\xi_3)$ and $\textbf{Q}=(\xi_1^2,\xi_2\xi_3)$. We assume $\textbf{Q}= (\xi_1^2,\xi_2^2+\xi_1\xi_3)$ firstly.  
Let $\xi=(\xi_1,\xi_2,\xi_3) \in \tau_1$ and $\eta=(\eta_1,\eta_2,\eta_3) \in \tau_2$. If
\begin{equation}\label{s6eq0}
	\Big|\textbf{n}^1_{\xi} \wedge \textbf{n}^2_{\xi}\wedge \textbf{n}^1_{\eta} \wedge \textbf{n}^2_{\eta} 	\Big|=0,
\end{equation}
by some calculations, it implies that $ \xi_1=\eta_1$.
Thus $\tau \in \mathcal{S}(B)$ is contained in one box with dimensions $K^{-1} \times 1 \times 1$. When we restrict $B_{K^2}$ in the physical space and $\mathcal{S}(B)$ in the frequency space, by locally constant property, $E^{\textbf{Q}}f$ behaves the same as $E^{\textbf{Q}_1}f$ with $\textbf{Q}_1=(0,\xi_2^2+\xi_1\xi_3)$. Then we use Theorem \ref{dec 3} on $\xi_3$-direction, one gets 
$$   \Big\|\sum_{\tau \in \mathcal{S}(B)}E^{\textbf{Q}}f_\tau\Big\|_{L^p(B)} \lesssim K^{\frac{1}{2}-\frac{1}{p}}\Big(   \sum_J \|E^{\textbf{Q}}f_J\|_{L^p(w_B)}^2\Big)^{\frac{1}{2}}, \quad \quad 2\leq p \leq \infty,  $$
where each $J$ denotes a box with dimensions $K^{-1} \times 1 \times K^{-1}$, and $w_B$ is given by (\ref{sig1}). When we restrict $B_{K^2}$ in the physical space and $J$ in the frequency space, by locally constant property, $E^{\textbf{Q}_1}f$ behaves the same as $E^{\textbf{Q}_2}f$ with $ \textbf{Q}_2=(0,\xi_2^2)$. Furtherly, we use Theorem \ref{dec 1} on the $\xi_2$-direction, 
$$   \|E^{\textbf{Q}}f_J\|_{L^p(w_B)} \lesssim K^{\beta(p)+\epsilon}\Big( \sum_{\tau \subset J} \|E^{\textbf{Q}}f_\tau\|_{L^p(w_B)}^2  \Big)^{\frac{1}{2}},   $$ 
where $\beta(p)=0$ when $2\leq p\leq 6$ and $\beta(p)=1/2-3/p$ when $6 \leq p\leq \infty$. 
Combining above estimates, we conclude
\begin{equation}\label{abc}
	\Big\|\sum_{\tau \in \mathcal{S}(B)}E^{\textbf{Q}}f_\tau\Big\|_{L^p(B)} \lesssim K^{\frac{1}{2}-\frac{1}{p}+\beta(p)+\epsilon} \Big( \sum_\tau \|E^{\textbf{Q}}f_\tau\|_{L^p(w_B)}^2  \Big)^{\frac{1}{2}}.
\end{equation}
We use (\ref{abc}) to run the argument of the narrow case as in Proposition \ref{sec2 le1} until the end. In fact, the proof becomes easier since it only has one scale $K$. Finally, we obtain  (\ref{aim}) with
$$  w(p) \geq \min_{p \geq 4 }  \max \Big\{   \frac{\alpha}{4},~   \frac{\beta(p)}{2}+\frac{\alpha}{4}+\frac{4-\alpha}{2p}   \Big\}=\min\Big\{\frac{\alpha+1}{4},\frac{\alpha+2}{6} \Big\}.   $$  
Then we assume $\textbf{Q}=(\xi_1^2,\xi_2\xi_3)$. By similar calculations, we have $\tau \in \mathcal{S}(B)$ is contained in a box with dimensions $K^{-1} \times 1 \times 1$. When we restrict $B_{K^2}$ in the physical space and $\mathcal{S}(B)$ in the frequency space,
by locally constant property, $E^{\textbf{Q}}f$ behaves the same as $E^{\textbf{Q}_3}f$ with  $ \textbf{Q}_3=(0,\xi_2\xi_3)$. We use Theorem \ref{dec 2} on the $\xi_2,\xi_3$-directions to get 
$$   \Big\|\sum_{\tau \in \mathcal{S}(B)}E^{\textbf{Q}}f_\tau\Big\|_{L^p(B)} \lesssim K^{\beta(p)+\epsilon} \Big( \sum_\tau \|E^{\textbf{Q}}f_\tau\|_{L^p(w_B)}^2  \Big)^{\frac{1}{2}},  $$
where $\beta(p)=1/2-1/p$ when $2\leq p\leq 6$ and $\beta(p)=1-4/p$ when $6 \leq p\leq \infty$.
Repeating the argument of the narrow case as in Proposition \ref{sec2 le1} until the end, we obtain  (\ref{aim}) with
$$ w(p)\geq  \min_{ p\geq 4 }  \max \Big\{   \frac{\alpha}{4},~  \frac{\beta(p)}{2}+\frac{\alpha-1}{4} +\frac{5-\alpha}{2p} \Big\}=\min\Big\{\frac{\alpha+1}{4},\frac{\alpha+2}{6} \Big\}.  $$ 
The argument of the items (b) and (c) in Theorem \ref{th1} is similar, here we omit it.

Now we study $\textbf{Q}$ satisfying (d) in Theorem \ref{th1}. It suffices to study  $\textbf{Q}=(\xi_1\xi_2,\xi_2^2+\xi_1 \xi_3)$ and $\textbf{Q}=(\xi_1\xi_2,\xi_1^2\pm\xi_3^2)$. We first consider the case  $\textbf{Q}=(\xi_1\xi_2,\xi_2^2+\xi_1 \xi_3)$. Let $\xi=(\xi_1,\xi_2,\xi_3) \in \tau_1$ and $\eta=(\eta_1,\eta_2,\eta_3) \in \tau_2$. If
\begin{equation}
	\Big|\textbf{n}^1_{\xi} \wedge \textbf{n}^2_{\xi}\wedge \textbf{n}^1_{\eta} \wedge \textbf{n}^2_{\eta} 	\Big|=0,
\end{equation}
by some calculations, it implies that 
$$ \xi_1=\eta_1, ~  \xi_2=\eta_2.$$
Therefore $\tau \in \mathcal{S}(B)$ is contained in a box with dimensions $K^{-1} \times K^{-1} \times 1$. When we restrict $B_{K^2}$ in the physical space and $\mathcal{S}(B)$ in the frequency space,
by locally constant property, $E^{\textbf{Q}}f$ behaves the same as $E^{\textbf{Q}_4}f$ with $\textbf{Q}_4 =( 0,\xi_1 \xi_3)$.  We use Theorem \ref{dec 3} on the $\xi_3$-direction to obtain 
\begin{equation}\label{s6eq5}
	\Big\|\sum_{\tau \in \mathcal{S}(B)}E^{\textbf{Q}}f_\tau\Big\|_{L^p(B)}\lesssim K^{\frac{1}{2}-\frac{1}{p}} \Big( \sum_\tau \|E^{\textbf{Q}}f_\tau\|_{L^p(w_B)}^2  \Big)^{\frac{1}{2}}, \quad \quad 2 \leq p \leq \infty. 
\end{equation}
Repeating the argument of the narrow case as in Proposition \ref{sec2 le1} until the end, we obtain (\ref{aim}) with
$$ w(p) \geq \min_{p \geq 4}  \max  \Big\{   \frac{\alpha}{4},~   \frac{\alpha}{4}+\frac{4-\alpha}{2p}   \Big\} =\frac{\alpha}{4}.$$  
Then we consider the case $\textbf{Q}=(\xi_1\xi_2,\xi_1^2\pm\xi_3^2)$. By testing (\ref{s6eq0}), we see that $\tau \in \mathcal{S}(B)$ is contained in a box with dimensions $K^{-1} \times K^{-1} \times 1$ or $K^{-1} \times 1\times K^{-1}$. In either case, we use Theorem \ref{dec 3} to obtain
\begin{equation}\label{s6eq7}
	\Big\|\sum_{\tau \in \mathcal{S}(B)}E^{\textbf{Q}}f_\tau\Big\|_{L^p(B)}\lesssim  K^{\frac{1}{2}-\frac{1}{p}+\epsilon}  \Big( \sum_\tau \|E^{\textbf{Q}}f_\tau\|_{L^p(w_B)}^2  \Big)^{\frac{1}{2}}, \quad \quad 2 \leq p \leq \infty. 
\end{equation}
Repeating the argument of the narrow case as in Proposition \ref{sec2 le1} until the end, we obtain (\ref{aim}) with  
$$ w(p) \geq \min_{p \geq 4}  \max  \Big\{   \frac{\alpha}{4},~   \frac{\alpha}{4}+\frac{4-\alpha}{2p}   \Big\} =\frac{\alpha}{4}.$$

Next we consider $\textbf{Q}$ satisfying (e) in Theorem \ref{th1}. \textbf{Claim}: the number of $\tau$'s in $\mathcal{S}(B)$ is bounded by $O(K)$.
Unlike the cases (a)-(d), we need more analysis to deal with the case (e) due to the absence of precise representations of $\textbf{Q}$. Let $\textbf{Q}=(Q_1,Q_2)$, $\xi=(\xi_1,\xi_2,\xi_3) \in \tau_1$ and $\eta=(\eta_1,\eta_2,\eta_3) \in \tau_2$. Note that the narrow condition
\begin{equation*}
	\Big|\textbf{n}^1_{\xi} \wedge \textbf{n}^2_{\xi}\wedge \textbf{n}^1_{\eta} \wedge \textbf{n}^2_{\eta} \Big|=0
\end{equation*}
 is equivalent to that $\textbf{n}^1_{\xi},~ ~\textbf{n}^2_{\xi},~\textbf{n}^1_{\eta},~\textbf{n}^2_{\eta}$ are linearly dependent, i.e. there exist $k_1$, $k_2$, $l_1$ and $l_2$ satisfying $k_1^2+k_2^2+l_1^2+l_2^2\neq 0$ such that
 $$   k_1\textbf{n}^1_{\xi} +k_2   \textbf{n}^2_{\xi} +l_1\textbf{n}^1_{\eta} +l_2\textbf{n}^2_{\eta} =0.  $$ 
This condition can be rewritten as:
\begin{equation*}
	\left(
	\begin{array}{cc}
		\partial_1 Q_1(\xi-\eta) & \partial_1 Q_2(\xi-\eta) \\
		\partial_2 Q_1(\xi-\eta) & \partial_2 Q_2(\xi-\eta) \\
			\partial_3 Q_1(\xi-\eta) & \partial_3 Q_2(\xi-\eta)
	\end{array}
	\right)
		\left(
	\begin{array}{c}
		k_1 \\
        k_2
	\end{array}
	\right)=0
\end{equation*}
has non-zero solutions on $(k_1, k_2)$. It follows that
\begin{equation}\label{s6bc1}
	\text{rank}\left(
	\begin{array}{cc}
		\partial_1 Q_1(\gamma) & \partial_1 Q_2(\gamma) \\
		\partial_2 Q_1(\gamma) & \partial_2 Q_2(\gamma) \\
		\partial_3 Q_1(\gamma) & \partial_3 Q_2(\gamma)
	\end{array}
	\right)\leq 1,
\end{equation}
where $\gamma=\xi-\eta$. 
If we set $Q_1(\gamma)= \gamma A \gamma^{T}$ and $Q_2(\gamma)= \gamma B \gamma^{T}$ where $A=(a_{ij})_{3\times 3}$ and $B=(b_{ij})_{3\times 3}$ are two real symmetric matrices, then it follows from (\ref{s6bc1}) that
\begin{equation}\label{s6bc2}
(a_{11} \gamma_1 +a_{12} \gamma_2 +a_{13} \gamma_3) (b_{21} \gamma_1 +b_{22} \gamma_2 +b_{23} \gamma_3)=(b_{11} \gamma_1 +b_{12} \gamma_2 +b_{13} \gamma_3)(a_{21} \gamma_1 +a_{22} \gamma_2 +a_{23} \gamma_3),
\end{equation}
and 
\begin{equation}\label{s6bc3}
(a_{11} \gamma_1 +a_{12} \gamma_2 +a_{13} \gamma_3) (b_{31} \gamma_1 +b_{32} \gamma_2 +b_{33} \gamma_3)=(b_{11} \gamma_1 +b_{12} \gamma_2 +b_{13} \gamma_3)(a_{31} \gamma_1 +a_{32} \gamma_2 +a_{33} \gamma_3).
\end{equation}
Denote $Z =\{\gamma: \gamma \text{~satisfies~} (\ref{s6bc2}) \text{~and~} (\ref{s6bc3})\}. $

We further do some reductions on $Q_1$ and $Q_2$ according to the condition of (e) in Theorem \ref{th1}. We pick $M\in \mathbb{R}^{3\times 3}$ and $M'\in \mathbb{R}^{1\times 2}$ such that the equality in (\ref{main sign}) is achieved with $d'=3,n'=1$. By linear transformations, we may suppose that $M=I_{3\times 3}$ and $M'=(1,0)$. Then $\mathfrak{d}_{3,1}(\textbf{Q})=2$ implies that $Q_1$ depends on 2 variables. We use linear transformations again to diagonalize $Q_1$, then
$$ (Q_1(\gamma),Q_2(\gamma))\equiv (\gamma_1^2\pm\gamma_2^2, b_{11}\gamma_1^2 +b_{22}\gamma_2^2+b_{33}\gamma_3^2+2b_{12}\gamma_1\gamma_2+2b_{13}\gamma_1\gamma_3+2b_{23}\gamma_2\gamma_3).     $$
It suffices to consider the case $Q_1(\gamma)=\gamma_1^2+\gamma_2^2$, since the argument of another case is similar. Moreover, $\mathfrak{d}_{3,2}(\textbf{Q})=3$ implies that $b_{13}$, $b_{23}$ and $b_{33}$ are not zero simultaneously. 
Firstly, we assume $b_{33}\neq 0$. By a change of variables in $\gamma_3$, we obtain
 $$ (Q_1(\gamma),Q_2(\gamma))\equiv (\gamma_1^2+\gamma_2^2, b_{11}\gamma_1^2 +b_{22}\gamma_2^2+2b_{12}\gamma_1\gamma_2+\gamma_3^2).     $$
In this setting, (\ref{s6bc2}) and (\ref{s6bc3}) become
\begin{equation}\label{s6bc4}
Z_1:~b_{12} (\gamma_1^2-\gamma_2^2)+(b_{22}-b_{11})\gamma_1 \gamma_2=0,
\end{equation}
and 
\begin{equation}\label{s6bc5}
Z_2:~	\gamma_1 \gamma_3=0.
\end{equation}
When we restrict $\gamma_1=0$, $Z_1$ becomes $b_{12}\gamma^2_2=0$, which is a line $\gamma_2=0$ as long as $b_{12}\neq 0$. In fact, we can assume that $b_{12}\neq0$. Otherwise, we have
  $$ (Q_1(\gamma),Q_2(\gamma))\equiv (\gamma_1^2+\gamma_2^2, b_{11}\gamma_1^2 +b_{22}\gamma_2^2+\gamma_3^2).     $$
If $b_{11}=b_{22}$, then adding a multiple of $Q_1$ to $Q_2$, one gets
  $$ (Q_1(\gamma),Q_2(\gamma))\equiv (\gamma_1^2+\gamma_2^2, \gamma_3^2),     $$
which is just the item (b) in Theorem \ref{th1}. If $b_{11}\neq b_{22}$, then adding a multiple of $Q_1$ to $Q_2$, one gets
  $$ (Q_1(\gamma),Q_2(\gamma))\equiv (\gamma_1^2+\gamma_2^2, \gamma_2^2+\gamma_3^2).     $$
Since such $\textbf{Q}$ possesses precise representation, we can verify (\ref{s6eq0}) directly and find the claim holds. When we restrict $\gamma_3=0$, $Z_1$ is a polynomial curve since we have excluded the case $b_{12}= 0$. Therefore we conclude
\begin{equation}\label{s6bc6}
Z=\{  \gamma_1=\gamma_2=0  \}   \cup  \{  b_{12} (\gamma_1^2-\gamma_2^2)+(b_{22}-b_{11})\gamma_1 \gamma_2=\gamma_3=0 \},
\end{equation} 
which is an algebraic variety with dimension 1. Next we assume $b_{33}=0$, so $b_{13}$ and $b_{23}$ are not zero simultaneously. Without loss of generality, we assume that $b_{13}\neq 0$, then
 $$ (Q_1(\gamma),Q_2(\gamma))\equiv (\gamma_1^2+\gamma_2^2, b_{11}\gamma_1^2 +b_{22}\gamma_2^2+2b_{12}\gamma_1\gamma_2+2\gamma_1 \gamma_3 +2 b_{23}\gamma_2 \gamma_3).     $$
In this setting, (\ref{s6bc2}) and (\ref{s6bc3}) become
 \begin{equation}\label{s6bc7}
 	Z_3:~b_{12} (\gamma_1^2-\gamma_2^2)+(b_{22}-b_{11})\gamma_1 \gamma_2+b_{23}\gamma_1\gamma_3-\gamma_2\gamma_3=0,
 \end{equation}
 and 
 \begin{equation}\label{s6bc8}
 	Z_4:~	\gamma_1 (\gamma_1+b_{23}\gamma_2)=0.
 \end{equation}
When we restrict $\gamma_1=0$, $Z_3$ becomes $\gamma_2 (b_{12}\gamma_2+\gamma_3)=0$, which a polynomial curve. When we restrict $\gamma_1+b_{23}\gamma_2=0$, $Z_3$ becomes
 $$ (b_{12} b_{23}^2-b_{12}+b_{22}b_{23}+b_{11}b_{23})\gamma^2_2-(1+b_{23}^2)\gamma_2\gamma_3=0.   $$
Note $1+b_{23}^2 \neq 0$, it is a polynomial curve. And then
\begin{align*}
&Z=\{  \gamma_1= \gamma_2 (b_{12}\gamma_2+\gamma_3)=0 \}   \\
 & \quad \quad \quad  \cup  \{  (b_{12} b_{23}^2-b_{12}+b_{22}b_{23}+b_{11}b_{23})\gamma^2_2-(1+b_{23}^2)\gamma_2\gamma_3=\gamma_1+b_{23}\gamma_2=0 \},
\end{align*}
which is also an algebraic variety with dimension 1. Thus we can use the Wongkew theorem \cite{W93} to obtain
$$   \Big|N_{K^{-1}}{Z} \cap [-1,1]^3 \Big|  \lesssim K^{-2}.  $$
This fact implies the claim immediately. Finally, we use Theorem \ref{dec 3} to obtain 
\begin{equation}
	\Big\|\sum_{\tau \in \mathcal{S}(B)}E^{\textbf{Q}}f_\tau\Big\|_{L^p(B)}\lesssim K^{\frac{1}{2}-\frac{1}{p}} \Big( \sum_\tau \|E^{\textbf{Q}}f_\tau\|_{L^p(w_B)}^2  \Big)^{\frac{1}{2}}, \quad \quad 2 \leq p \leq \infty. 
\end{equation}
Repeating the argument of the narrow case as in Proposition \ref{sec2 le1} until the end, we obtain (\ref{aim}) with
$$ w(p) \geq \min_{p \geq 4}  \max  \Big\{   \frac{\alpha}{4},~   \frac{\alpha}{4}+\frac{4-\alpha}{2p}   \Big\} =\frac{\alpha}{4}.$$

\qed

\begin{remark}
We explain why the opposite of (\ref{s6e1}) is just (\ref{s6eq0}). Since the scale of $\tau$ is $K^{-1}$, if (\ref{s6eq0}) fails, there exists one constant $C>0$ such that
\begin{equation*}
	\Big| \mathbf{n}^1_{\xi} \wedge \mathbf{n}^2_{\xi}\wedge \mathbf{n}^1_{\eta} \wedge \mathbf{n}^2_{\eta}  	\Big| \gtrsim K^{-C}
\end{equation*} 
for any $\xi \in \tau_1$, $\eta \in \tau_2$. For every quadratic form $\textbf{Q}$, we can determine suitable $C$ by some simple calculations. For example, for $\textbf{Q}= (\xi_1^2,\xi_2^2+\xi_1\xi_3)$, there is
$$   	\Big| \mathbf{n}^1_{\xi} \wedge \mathbf{n}^2_{\xi}\wedge \mathbf{n}^1_{\eta} \wedge \mathbf{n}^2_{\eta}  	\Big| \sim |\xi_1-\eta_1|,  $$
 and we can take $C=1$.  
\end{remark}

Through the similar argument as in Section 3, we obtain the following  weighted restriction estimates. 
\begin{proposition}\label{cor3}
Let $\alpha \in (0,5]$, and $\mu$ be an $\alpha$-dimensional measure in $\mathbb{R}^5$. Then
	\begin{itemize}
\item[(1)] For $\textbf{Q}$ satisfying ${\rm (a)}$ in Theorem \ref{th1}, we have
\begin{equation}
	\|E^{\textbf{Q}}f\|_{L^2(d\mu_R)}\lesssim R^{\min\{   \frac{\alpha+1}{4},\frac{\alpha+2}{6} \}  +\epsilon}\|f\|_{L^2}.   
\end{equation}
\item[(2)] For $\textbf{Q}$ satisfying ${\rm (b)}$ in Theorem \ref{th1}, we have
\begin{equation}
	\|E^{\textbf{Q}}f\|_{L^2(d\mu_R)}\lesssim R^{\max\left\{ \frac{\alpha}{4},  \min\{    \frac{\alpha+1}{4}, \frac{\alpha+3}{8} \}    \right\}+\epsilon}\|f\|_{L^2}.  
\end{equation}
\item[(3)] For $\textbf{Q}$ satisfying ${\rm (c)}$ in Theorem \ref{th1}, we have
\begin{equation}
	\|E^{\textbf{Q}}f\|_{L^2(d\mu_R)}\lesssim R^{ \frac{\alpha+1}{4} +\epsilon}\|f\|_{L^2}.  
\end{equation}
		\item[(4)] For $\textbf{Q}$ satisfying ${\rm (d)}$ and ${\rm (e)}$ in Theorem \ref{th1}, we have
		\begin{equation}\label{special}
			\|E^{\textbf{Q}}f\|_{L^{2}(d\mu_R)}\lesssim  R^{\frac{\alpha}{4}+\epsilon} \|f\|_{L^2}.
		\end{equation}
	\end{itemize}
\end{proposition}

\section{Stein-Tomas type's inequalities}

In the previous sections, we have adopted several methods to obtain a few  weighted restriction estimates for quadratic forms. We recall the notation $s(\alpha,S_{\textbf{Q}})$, which denotes the optimal constant $s$ such that the inequality
\begin{equation*}
	\|E^{\textbf{Q}}f\|_{L^2(d\mu_R)}\lesssim R^{s}\|f\|_{L^2}
\end{equation*}
holds for every $f \in L^2(\mathbb{R}^d)$, every $\alpha$-dimensional measure $\mu\in\mathcal{M}(B^{d+n}(0,1))$, and  any $R \geq 1$. In general, we expect $s(\alpha,S_{\textbf{Q}}) \rightarrow 0$ when $\alpha \rightarrow 0$. However, the results obtained in the previous sections are not very nice except (\ref{special}) for small $\alpha$. To obtain better corresponding bounds for small $\alpha$, 
we first prove the almost sharp Stein-Tomas type's inequalities for the quadratic forms in Theorem \ref{th1}.

\begin{theorem}[Stein-Tomas type's inequalities for surfaces of co-dimension 2 in $\mathbb{R}^5$]
Let $d=3$ and $n=2$. For any $R \geq 1$, the following holds:
	\begin{itemize}
		\item[(1)] For $\textbf{Q}$ satisfying ${\rm (a)}$, ${\rm (b)}$ and ${\rm (c)}$ in Theorem \ref{th1}, we have
		\begin{equation}
			\|E^{\textbf{Q}}f\|_{L^6(B_R)}\lesssim R^{\epsilon}\|f\|_{L^2}.  
		\end{equation}
		\item[(2)] For $\textbf{Q}$ satisfying ${\rm (d)}$ in Theorem \ref{th1}, we have
		\begin{equation}\label{special2}
			\|E^{\textbf{Q}}f\|_{L^{\frac{16}{3}}(B_R)}\lesssim  R^{\epsilon} \|f\|_{L^2}, \quad   \textbf{Q}=(\xi_1 \xi_2,\xi_2^2+\xi_1 \xi_3),
		\end{equation}
		and
		\begin{equation}
			\|E^{\textbf{Q}}f\|_{L^5(B_R)}\lesssim R^{\epsilon} \|f\|_{L^2}, \quad   \textbf{Q}=(\xi_1 \xi_2,\xi_1^2\pm\xi_3^2).
		\end{equation}
		\item[(3)] For $\textbf{Q}$ satisfying ${\rm (e)}$ in Theorem \ref{th1}, we have
		\begin{equation}
			\|E^{\textbf{Q}}f\|_{L^{\frac{14}{3}}(B_R)}\lesssim \|f\|_{L^2}. 
		\end{equation}
	\end{itemize}
Moreover, the estimates of ${\rm (1)}$ and ${\rm (2)}$ are sharp up to the endpoints, and the estimate of ${\rm (3)}$ is sharp.
\end{theorem}

\noindent \textit{Proof.} The estimates of (1) and (3) are implied in the argument of \cite{GO} and \cite{M96}, respectively. We only need to consider (2).
If  $\textbf{Q}=(\xi_1 \xi_2,\xi_2^2+\xi_1 \xi_3)$, we run the argument of \cite[Section 6]{GO}, then (\ref{special2}) follows. Therefore, it suffices to deal with $\textbf{Q}=(\xi_1 \xi_2,\xi_1^2\pm\xi_3^2)$. Now we find what $p\geq 2$ such that the following estimate holds:
\begin{equation}\label{7aim}
\|E^{\textbf{Q}}f\|_{L^p(B_R)} \lesssim R^\epsilon \|f\|_{L^2}.
\end{equation}  

We run the broad-narrow analysis once again. Divide $B^3(0,1)$ to $K^{-1}$-balls $\tau$, and write $f=\sum_\tau f_\tau$, where $f_\tau=f\chi_\tau$. For each $K^2$-ball $B$ in $B_R$, define its significant set as
$$    \mathcal{S}(B) :=\Big\{ \tau: \| E^{\textbf{Q}} f _\tau\|_{L^p(B)} \geq \frac{1}{100\# \tau} \|  E^{\textbf{Q}}f  \|_{L^p(B)} \Big\}.  $$
Then
$$    \Big\|  \sum_{\tau \in \mathcal{S}(B)} E^{\textbf{Q}}f_\tau \Big\|_{L^p(B)} \sim \|E^{\textbf{Q}}f\|_{L^p(B)}.  $$
We say $B$ is broad if there exist $\tau_1,\tau_2 \in \mathcal{S}(B)$ such that for any $\xi \in \tau_1$, $\eta \in \tau_2$
\begin{equation}\label{s7e1}
		\Big| \textbf{n}^1_{\xi} \wedge \textbf{n}^2_{\xi}\wedge \textbf{n}^1_{\eta} \wedge \textbf{n}^2_{\eta}  	\Big| \gtrsim K^{-O(1)}.
\end{equation}
Otherwise, we say $B$ is narrow. We denote the union of broad balls $B $ in $ B_R$ by $Y_{\text{broad}}$ and the union of narrow balls $B $ in $  B_R$ by $Y_{\text{narrow}}$. 

Firstly, we consider the contribution from $Y_{\text{broad}}$. For each broad ball $B$, there exist $\tau_1,\tau_2 \in \mathcal{S}(B) $ satisfying the transverse condition (\ref{s7e1}) such that
$$   \int_B  |E^{\textbf{Q}}f|^p \lesssim   K^{3p} \prod_{j=1}^{2} \Big(  \int_B |E^{\textbf{Q}}f_{\tau_j}|^p  \Big)^{\frac{1}{2}}.  $$
Repeating the argument as in Section 3, together with locally constant property, one concludes
\begin{align*}
\int_B  |E^{\textbf{Q}}f|^p & \lesssim K^{O(1)} \int_{B(x_B,2)} \prod_{j=1}^{2} |E^{\textbf{Q}}f_{\tau_j,v_j}|^{\frac{p}{2}}  \\
& \leq  K^{O(1)} \sum_{\tau_1,\tau_2~\text{trans}}  \sum_{v_1,v_2} \int_{B(x_B,2)} \prod_{j=1}^{2} |E^{\textbf{Q}}f_{\tau_j,v_j}|^{\frac{p}{2}},
\end{align*}
where $v_j \in B(0,K^2) \cap \mathbb{Z}^{5}$.
Summing all $B $ in $Y_{\text{broad}}$, together with Theorem \ref{klin res} for $p \geq 4$, one has
\begin{align*}
	\int_{Y_{\text{broad}} }  |E^{\textbf{Q}}f|^p  \lesssim K^{O(1)}\sum_{\tau_1,\tau_2~\text{trans}}  \sum_{v_1,v_2} \int_{B_R} \prod_{j=1}^{2} |E^{\textbf{Q}}f_{\tau_j,v_j}|^{\frac{p}{2}}   \lesssim  R^{p\epsilon}\|f\|_{L^2}^p.
\end{align*}

Next we consider the contribution from $Y_{\text{narrow}}$. For each narrow ball $B$, we have
$$  \|E^{\textbf{Q}}f\|_{L^p(B)} \lesssim \Big\|  \sum_{\tau \in \mathcal{S}(B)}E^{\textbf{Q}}f_\tau \Big\|_{L^p(B)}.   $$
We now test the opposite of (\ref{s7e1}). For two points $\xi=(\xi_1,\xi_2,\xi_3) \in \tau_1$ and $\eta=(\eta_1,\eta_2,\eta_3) \in \tau_2$, if 
$$	\Big|\textbf{n}^1_{\xi} \wedge \textbf{n}^2_{\xi}\wedge \textbf{n}^1_{\eta} \wedge \textbf{n}^2_{\eta}	\Big|=0,$$ 
by some calculations, it implies that
$$  \xi_1=\eta_1, ~\xi_2=\eta_2 \quad \text{or} \quad  \xi_1=\eta_1, ~\xi_3=\eta_3 . $$ 
The above relation shows that $\tau \in \mathcal{S}(B)$ is contained in a box with dimensions $K^{-1} \times K^{-1} \times 1$ or $K^{-1} \times 1\times K^{-1}$.
 We denote the union of narrow balls $B$ in the first case by $Y_{1}$ and the union of narrow balls $B$ in the second case by $Y_{2}$.

For each narrow ball $B $ in $Y_1$, we have $\tau \in \mathcal{S}(B)$ is contained in a box of dimensions $K^{-1} \times K^{-1} \times 1$. When we restrict $B_{K^2}$ in the physical space and $\mathcal{S}(B)$ in the frequency space,
by locally constant property, $E^{\textbf{Q}}f$ behaves the same as $E^{\textbf{Q}_5}f$ with $\textbf{Q}_5 =( 0,\xi_3^2)$. Using Theorem \ref{dec 1} on the $\xi_3$-direction, for $2\leq p \leq 6$, one has
$$  \|E^{\textbf{Q}}f\|_{L^p(B)} \lesssim K^{\frac{\epsilon}{2}}\Big(\sum_{\tau \in \mathcal{S}(B)}\left\|  E^{\textbf{Q}}f_\tau \right\|^2_{L^p(w_B)}\Big)^{\frac{1}{2}}\leq  K^{\frac{\epsilon}{2}} \Big(\sum_{\tau }\left\|  E^{\textbf{Q}}f_\tau \right\|^2_{L^p(w_B)}\Big)^{\frac{1}{2}}.  $$
Summing all $B $ in $ Y_1$ by the parallel decoupling inequality, we get 
$$    \|E^{\textbf{Q}}f\|_{L^p(Y_1)} \lesssim K^{\frac{\epsilon}{2}}\Big(\sum_{\tau}\left\|  E^{\textbf{Q}}f_\tau \right\|^2_{L^p(B_R)}\Big)^{\frac{1}{2}}.   $$
As for the term $\|  E^{\textbf{Q}}f_\tau \|_{L^p(B_R)}$, using the change of coordinates
$$   \xi   \rightarrow c(\tau)+K^{-1}\xi, $$
we conclude
\begin{equation*}
	\|E^{\textbf{Q}} f _\tau\|_{L^p(B_R)}=K^{\frac{7}{p}-\frac{3}{2}} \|E^{\textbf{Q}}g\|_{L^p(\tilde{Y})},
\end{equation*}
for some $g$ with supp$g\subset B^3(0,1)$ and $\|g\|_{L^2}=\|f_\tau\|_{L^2}$, where $\tilde{Y}$ is contained in a ball with radius $R/K$. Applying inductive hypothesis of (\ref{7aim}) at scale $R/K$, we conclude
$$  	\|E^{\textbf{Q}} f \|_{L^p(Y_1)} \lesssim K^{\frac{7}{p}-\frac{3}{2}-\frac{\epsilon}{2}}R^\epsilon \|f\|_{L^2}. $$
Therefore, the induction closes of $Y_1$ when $p \geq 14/3$.

For each narrow ball $B$ in $Y_2$, we have $\tau \in \mathcal{S}(B)$ is contained in a box of dimensions $K^{-1} \times 1 \times K^{-1}$. We denote such boxes by $J$. Then
$$    \|E^{\textbf{Q}}f\|_{L^p(B)} \leq \Big(  \sum_J \|E^{\textbf{Q}}f_J\|_{L^p(B)}^2    \Big)^{\frac{1}{2}}.      $$
Summing all $B$ in $ Y_2$ by the parallel decoupling inequality, we get 
$$    \|E^{\textbf{Q}}f\|_{L^p(Y_2)} \leq \Big(  \sum_J \|E^{\textbf{Q}}f_J\|_{L^p(B_R)}^2    \Big)^{\frac{1}{2}}.      $$
Set $c(\tau)=(c_1(\tau),c_2(\tau),c_3(\tau))$. Using the change of coordinates
$$   \xi_1   \rightarrow c_1(\tau)+K^{-1}\xi_1,~ \xi_2   \rightarrow c_2(\tau)+\xi_2,~ \xi_3   \rightarrow c_3(\tau)+K^{-1}\xi_3, $$
we conclude
\begin{equation*}
	\|E^{\textbf{Q}} f _J\|_{L^p(B_R)}=K^{\frac{5}{p}-1} \|E^{\textbf{Q}}g\|_{L^p(\tilde{Y})},
\end{equation*}
for some $g$ with supp$g\subset B^3(0,1)$ and $\|g\|_{L^2}=\|f_J\|_{L^2}$, where $\tilde{Y}$ is contained in a box with dimensions $R/K \times R \times R/K$. To apply the induction hypothesis, we break $\tilde{Y}$ into smaller balls of radius $R/2$, and obtain
$$  	\|E^{\textbf{Q}} f \|_{L^p(Y_2)} \lesssim K^{\frac{5}{p}-1+\epsilon}R^\epsilon \|f\|_{L^2}. $$
Therefore, the induction closes of $Y_2$ when $p >5$.

Combining all above constraints of $p$, we conclude (\ref{7aim}) holds for $5 <p \leq 6$. Then (\ref{7aim}) for $p=5$ is a direct result  by H\"older's inequality.

Finally, we point out the values of $p$ in (1)-(3) are sharp, where we can easily verify this fact from Table 1.

\qed

As a corollary of this result, we can obtain the relevant weighted restriction estimates by using H\"older's inequality.

\begin{corollary}\label{cor4}
Let $\alpha \in (0,5]$, and $\mu$ be an $\alpha$-dimensional measure in $\mathbb{R}^5$. Then
	\begin{itemize}
		\item[(1)] For $\textbf{Q}$ satisfying ${\rm (a)}$, ${\rm (b)}$ and ${\rm (c)}$ in Theorem \ref{th1}, we have
		\begin{equation}
			\|E^{\textbf{Q}}f\|_{L^2(d\mu_R)}\lesssim R^{\frac{\alpha}{3}+\epsilon}\|f\|_{L^2}.  
		\end{equation}
		\item[(2)] For $\textbf{Q}$ satisfying ${\rm (d)}$ in Theorem \ref{th1}, we have
		\begin{equation}
			\|E^{\textbf{Q}}f\|_{L^{2}(d\mu_R)}\lesssim  R^{\frac{5\alpha}{16}+\epsilon} \|f\|_{L^2}, \quad   \textbf{Q}=(\xi_1 \xi_2,\xi_2^2+\xi_1 \xi_3),
		\end{equation}
		and
		\begin{equation}
			\|E^{\textbf{Q}}f\|_{L^2(d\mu_R)}\lesssim R^{\frac{3\alpha}{10}+\epsilon} \|f\|_{L^2}, \quad   \textbf{Q}=(\xi_1 \xi_2,\xi_1^2\pm\xi_3^2).
		\end{equation}
		\item[(3)] For $\textbf{Q}$ satisfying ${\rm (e)}$ in Theorem \ref{th1}, we have
		\begin{equation}
			\|E^{\textbf{Q}}f\|_{L^{2}(d\mu_R)}\lesssim R^{\frac{2\alpha}{7}}\|f\|_{L^2}. 
		\end{equation}
	\end{itemize}
\end{corollary}

Finally, we have showed Theorem \ref{th3} by combining Section 2.1, Corollary \ref{cor1}, Corollary \ref{cor2}, Proposition \ref{cor3}, Corollary \ref{cor4} and the following Agmon-H\"ormander inequality:
\begin{equation}
		\|E^{\textbf{Q}}f\|_{L^2(B_R)} \lesssim R\|f\|_{L^2}.
\end{equation}

\appendix

\section{The refined Strichartz estimates for quadratic forms}

The refined Strichartz estimate is another important tool to establish the weighted restriction estimates. It was used to gain almost sharp pointwise convergence range of the solution of Schr\"odinger equation in $\mathbb{R}^3$ in \cite{DGL} for the first time. Now we extend their results to the cases of general quadratic forms, where the proof is essentially the same as in \cite{DGLZ}.
\begin{theorem}[Refined Strichartz estimates for quadratic forms]\label{rseqf}
Let $d,n \geq 1$ and $p \geq 2$. Let  $\textbf{Q}=(Q_1,...,Q_n)$ be a sequence of quadratic forms defined on $\mathbb{R}^d$.  Suppose that $f$ satisfies ${\rm supp}f \subset B^d(0,1)$. Suppose that $Q_1,Q_2,...$ are lattice $R^{1/2}$-cubes in $B^{d+n}(0,R)$, and 
	$$    \|E^{\textbf{Q}} f\|_{L^{p}(Q_j)} \text{~is~essentially~a ~dyadic~constant~in~}  j.   $$
	Suppose that these cubes are arranged in horizontal slabs of the form $\mathbb{R} \times ...\times \mathbb{R} \times [t_1,t_1+R^{1/2}] \times... \times [t_n,t_n+R^{1/2}]$, and each such strip contains $\sim \sigma$ cubes $Q_j$. Let $Y=\cup_j Q_j$, then for any $\epsilon>0$, we have
	\begin{equation}
		\left\|  E^{\textbf{Q}} f \right\|_{L^{p}(Y)}   \leq C_\epsilon R^{w(p)+\epsilon}\sigma^{\frac{1}{p}-\frac{1}{2}}\|f\|_{L^2}, 
	\end{equation}
	where $w(p)$ is given by 
	$$   w(p)=\frac{1}{2}   \Big[\Gamma_p^{d}(\textbf{Q})+\frac{d+2n}{p}-\frac{d}{2}    \Big]       .    $$
	In particular, if $\textbf{Q}$ is strongly non-degenerate, for $p=2(d+2n)/d$, we have
	\begin{equation}
		\left\|  E^{\textbf{Q}} f \right\|_{L^{\frac{2(d+2n)}{d}}(Y)}   \leq C_\epsilon R^{\epsilon}\sigma^{-\frac{n}{d+2n}}\|f\|_{L^2}. 
	\end{equation}
\end{theorem}

Corollary \ref{buchong2} can be also derived by Theorem \ref{rseqf} through the same argument as in \cite[Section 3]{DGLZ}. This phenomenon seems natural since the methods from Section 4 and this section both rely on the  original dimensional $\ell^2$ decoupling.

For completeness, we point out that, via the same technique as in \cite{DGL,DGLZ}, together with Theorem \ref{k linear kakeya}, we can obtain the following $k$-linear refined Strichartz estimates for quadratic forms.

\begin{theorem}[$k$-linear refined Strichartz estimates for quadratic forms]
Let $d,n \geq 1$ and $p \geq 2$. Let  $\textbf{Q}=(Q_1,...,Q_n)$ be a sequence of quadratic forms defined on $\mathbb{R}^d$.  Suppose that $f_i,i=1,...,k,$ satisfy ${\rm supp}f_i \subset B^d(0,1)$ and the $k$-transverse condition (\ref{k trans con}). Suppose that $Y=\cup_{j=1}^N Q_j$ are lattice $R^{1/2}$-cubes in $B^{d+n}(0,R)$, and 
	$$    \|E^{\textbf{Q}} f_i\|_{L^{p}(Q_j)} \text{~is~essentially~a ~dyadic~constant~in~}    j,  \quad i=1,...,k. $$
	Then for any $\epsilon>0$, we have
	\begin{equation}
		\left\| \prod_{i=1}^{k} |E^{\textbf{Q}} f_i|^{\frac{1}{k}} \right\|_{L^{p}(Y)}   \leq C_\epsilon R^{w(p)+\epsilon}N^{-\frac{k-1}{k}(\frac{1}{2}-\frac{1}{p})}\prod_{i=1}^{k}\|f_i\|^{\frac{1}{k}}_{L^2}, 
	\end{equation}
	where $w(p)$ is given by 
	$$   w(p)=\frac{1}{2}   \Big[\Gamma_p^{d}(\textbf{Q})+\frac{d+2n}{p}-\frac{d}{2}    \Big]       .    $$
		In particular, if $\textbf{Q}$ is strongly non-degenerate, for $p=2(d+2n)/d$, we have
	\begin{equation}
		\left\| \prod_{i=1}^{k} |E^{\textbf{Q}} f_i|^{\frac{1}{k}} \right\|_{L^{\frac{2(d+2n)}{d}}(Y)}   \leq C_\epsilon R^{\epsilon}N^{-\frac{(k-1)n}{k(d+2n)}}\prod_{i=1}^{k}\|f_i\|^{\frac{1}{k}}_{L^2}. 
	\end{equation}
\end{theorem}

\section{the proof of Theorem \ref{th1}}

In this appendix, we prove Theorem \ref{th1}. In fact, Guo and Oh \cite{GO} have classified quadratic forms of co-dimension 2 in $\mathbb{R}^5$ according to the (D)  condition and $({\rm R_{n'}})$ condition.  We say that $\textbf{Q}=(Q_1,Q_2)$ satisfies the (D) condition if
\begin{equation}\label{ap21}
\left|
	\begin{array}{ccc}
		\partial_1 Q_1(\xi) & \partial_2 Q_1(\xi) &\partial_3 Q_1(\xi) \\
		\partial_1 Q_2(\xi) & \partial_2 Q_2(\xi) &\partial_3 Q_2(\xi)\\
		w_1 & w_2 & w_3
	\end{array}
	\right|\not\equiv 0,
\end{equation}
for any non-zero vector $w=(w_1,w_2,w_3)\in \mathbb{R}^3$. We say that $\textbf{Q}=(Q_1,Q_2)$ satisfies the $({\rm R_{n'}})$ condition if
\begin{equation}\label{ap22}
\inf_{\text{hyperplane~} H\subset \mathbb{R}^3}  \sup_{\lambda_1,\lambda_2 \in \mathbb{R}}  \text{rank}\Big(  (\lambda_1Q_1+\lambda_2 Q_2)|_H  \Big)=n'.
\end{equation}

We only need to show the relations between the two conditions and the values of $\mathfrak{d}_{d',n'}(\textbf{Q})$. On the one hand, by Lemma 2.2 in \cite{GORYZK19}, combined with the items (a) and (b) in Theorem \ref{th1}, we easily see that
$$   \textbf{Q}  \text{~satisfies~ the~} ({\rm D}) \text{~condition~}\Longleftrightarrow ~ \mathfrak{d}_{3,1}(\textbf{Q})>1.  $$
On the other hand, by (1.14) in \cite{GOZZK}, we easily see that 
$$  \textbf{Q}  \text{~satisfies~ the~} ({\rm R_{n'}}) \text{~condition~}\Longleftrightarrow ~\mathfrak{d}_{2,2}(\textbf{Q})=n'. $$
Finally, we prove that $\mathfrak{d}_{3,1}(\textbf{Q})\leq 2$ for all $\textbf{Q}$. Fix a quadratic form $\textbf{Q}=(Q_1,Q_2)$ in 3 variables.
Let $a$, $b$ are two real numbers with $a^2+b^2\neq 0$, then
$$     aQ_1(\xi) +b Q_2(\xi) =   \xi M(a,b) \xi^{T},    $$ 
where $M(a,b)$ is a real symmetric matrix. We observe that $\text{det}M(a,b)$ is a homogeneous polynomial of degree at most 3. Without loss of generality, we can assume the degree of $\text{det}M(a,b)$ is 3. Otherwise, at least one of $Q_1$ and $Q_2$ is non-invertible, and  $\mathfrak{d}_{3,1}(\textbf{Q})\leq 2$ obviously.  By the fundamental theorem of algebra, we can write
$$  \text{det}M(a,b)=(\lambda a+\mu b) P(a,b). $$
Here $\lambda$, $\mu$ are two real numbers satisfying $\lambda^2+\mu^2\neq0$, and $P(a,b)$ is a homogeneous polynomial of degree 2. Thus there exist $a_0$ and $b_0$ with $a_0^2+b_0^2\neq 0$ such that $\text{det}M(a_0,b_0)=0$. Choosing one appropriate real invertible matrix $M \in  \mathbb{R}^{3\times 3}$ to diagonalize $M(a_0,b_0)$, we obtain
$$	\text{NV}(    (a_0 Q_1+b_0 Q_2)\circ M   ) \leq 2.$$
By the definition of $\mathfrak{d}_{d',n'}(\textbf{Q})$, we obtain $\mathfrak{d}_{3,1}(\textbf{Q})\leq 2$.

\vskip0.3cm

\subsection*{Acknowledgements} 
We thank the anonymous referee and the associated editor for their invaluable comments which helped to improve the paper. This work was supported by the National Key R${\&}$D Program of China (No. 2022YFA1005700), NSFC Grant No. 12371095  and No. 12371100.

\bibliographystyle{plain}
\bibliography{mybibfile}

\end{document}